\documentclass[12pt,a4paper]{article}
\usepackage{amsfonts}
\usepackage{mathrsfs}
\usepackage{amsmath}
\usepackage{amssymb}

\usepackage{bbm}
\usepackage{mathrsfs}
\usepackage{graphicx}

\usepackage{enumerate}
\usepackage{theorem}

\makeatletter
\def\tank#1{\protected@xdef\@thanks{\@thanks
        \protect\footnotetext[0]{#1}}}
\def\bigfoot{

    \@footnotetext}
\makeatother

\newcommand{\ea}{\end{array}}
\newtheorem{theorem}{Theorem}[section]
\newtheorem{proposition}{Proposition}[section]
\newtheorem{corollary}{Corollary}[section]
\newtheorem{lemma}{Lemma}[section]
\newtheorem{definition}{Definition}[section]
\newtheorem{Rem}{Remark}[section]

{\theorembodyfont{\rmfamily}
}

\newenvironment{proof}{Proof.}

\begin{document}
\title{\Large\bf 2D Stochastic Chemotaxis-Navier-Stokes System}
\footnotesize{\author{Jianliang Zhai $^{1,}$\thanks{zhaijl@ustc.edu.cn},\ \
 Tusheng Zhang$^{2}$\thanks{Tusheng.Zhang@manchester.ac.uk}\\
 {\em $^1$ School of Mathematical Sciences,}\\
 {\em University of Science and Technology of China,}\\
 {\em Hefei, 230026, China}\\
 {\em $^2$ School of Mathematics, University of Manchester,}\\
 {\em Oxford Road, Manchester, M13 9PL, UK}\\
}
\maketitle

\noindent{\bf Abstract:}
In this paper, we establish the existence and uniqueness of both mild(/variational) solutions and weak (in the sense of PDE) solutions of coupled system of 2D stochastic Chemotaxis-Navier-Stokes equations. The mild/variational solution is obtained through a fixed point argument in a purposely constructed Banach space. To get the weak solution we  first prove the existence of a martingale weak solution and then we show that the pathwise uniqueness holds for the martingale solution.

\vskip 0.4cm

\vskip 0.4cm


\tableofcontents

\section{Introduction}
\setcounter{equation}{0}
The purpose of this paper is to establish the existence and uniqueness of the solution of the coupled 2D stochastic Chemotaxis-Navier-Stokes system:
\begin{eqnarray}\label{eq system 00}
&& dn+u\cdot \nabla ndt=\delta \Delta ndt-\nabla\cdot(\chi(c)n\nabla c)dt,\nonumber\\
&& dc+u\cdot \nabla cdt=\mu \Delta cdt-k(c)ndt,\nonumber\\
&& du+(u\cdot \nabla)udt+\nabla Pdt=\nu \Delta udt-n\nabla\phi ~dt+\sigma(u)dW_t,\\
&& \nabla\cdot u=0,\ \ \ \ \ t>0,\ x\in\mathcal{O}.\nonumber
\end{eqnarray}
The system arises  in the modeling of bacterial suspensions in fluid drops and describes the spontaneous emergence of patterns in populations of oxygen-driven swimming bacteria.
Here, $\mathcal{O}\subset\mathbb{R}^2$ is a bounded convex domain with smooth boundary
$\partial \mathcal{O}$, which will be
 the  spatial domain where the moving cells and the fluid interact. The unknowns are $n=n(t,x):\mathbb{R}^+\times\mathcal{O}\rightarrow\mathbb{R}^+$, $c(t,x):\mathbb{R}^+\times\mathcal{O}\rightarrow\mathbb{R}^+$,
$u(t,x):\mathbb{R}^+\times\mathcal{O}\rightarrow\mathbb{R}^2$ and $P=P(t,x): \mathbb{R}^+\times\mathcal{O}\rightarrow\mathbb{R}$,
 which represent respectively the cell density, chemical concentration, velocity field and pressure of the fluid.  Positive constants $\delta,\mu,\nu$ are the corresponding diffusion coefficients for the cells, chemical and fluid. The gravitational potential $\phi=\phi(x)$, the chemotactic sensitivity $\chi(c)$ and the per-capita oxygen consumption rate $k(c)$ are supposed to be given sufficiently smooth functions.   $\{W_t,\ t\geq0\}$ is a cylindrical  Wiener process
representing the external random driving force.

System (\ref{eq system 00}) is considered  with the boundary conditions
\begin{eqnarray}\label{eq boundary condition 1}
\frac{\partial n}{\partial v}=\frac{\partial c}{\partial v}=0\text{  and  }u=0\ \ \text{   for }x\in\partial\mathcal{O}\text{ and }t>0,
\end{eqnarray}
and the initial conditions
\begin{eqnarray}\label{eq boundary condition 2}
n(0,x)=n_0(x),\ \ \ c(0,x)=c_0(x),\ \ \ u(0,x)=u_0(x),\ \ \ x\in\mathcal{O}.
\end{eqnarray}
\vskip 0.4cm

The deterministic models of system (\ref{eq system 00}) (i.e. $\sigma=0$) was proposed by Tuval \emph{et al.} in \cite{TCDWKG 2005}. In \cite{Xue Othmer 2009}, the authors suggest a wider variants to describe more complicated interaction neighborhood environment around cells. The well-posedness of the deterministic models of system (\ref{eq system 00}) (and its variants) is a highly non-trivial problem. In the past several years, the main focus of the existing literature is on the solvability of the system, see \cite{Cao 2016,CKL 2014,Duan Xiang 2014, Ishida 2015, KMS 2016, Li Li 2016, Liu-Lorz, Tao Winkler 2013, Winkler, Winkler 10, Winkler 2016, Winkler 2015 01, Zhang Zheng 2014} and reference therein. We like to mention a few of them which are relevant to our work.  In \cite{Lorz 2010}, local (in time) weak solutions (in the sense of PDE) were constructed in a bounded domain in $\mathbb{R}^d$, $d=2,3$ with no-flux boundary condition and in $\mathbb{R}^2$ for a special case. Based on some nice energy estimates, if the convective term $(u\cdot \nabla)u$ is neglected, global weak solutions were obtained in \cite{Duan Lorz Markowich 2010} provided the initial data or $\nabla \phi$ is small. Our work is motivated and influenced by the  recent papers \cite{Liu-Lorz} and \cite{Winkler}. In \cite{Liu-Lorz}, for the models in $\mathbb{R}^2$, Liu and Lorz developed some nice entropy estimates  to prove the global existence of weak solutions to the deterministic models of system (\ref{eq system 00}) for large initial data. In \cite{Winkler}, when $\mathcal{O}\subset\mathbb{R}^2$ is a bounded convex domain with smooth boundary
$\partial \mathcal{O}$, the author managed to establish the existence and uniqueness of  global strong (in the sense of PDE) solution of system (\ref{eq system 00}) without the restriction of the smallness of either the initial data or the coefficients. There are many other interesting results on this topic, we refer to the references mentioned above. Finally, we refer the reader to \cite{Winkler 2015 00,Winkler 2014 00} for the stabilization and convergence rate of solutions of the deterministic models of system (\ref{eq system 00}) and its variants.
\vskip 0.3cm

Taking into account the  random environment the bacteria are in  and the effect of  random external forces, it is natural to consider the coupled 2D stochastic Chemotaxis-Navier-Stokes system (\ref{eq system 00}). Adding the singular random noise to the system  changes the mathematical analysis significantly. In this paper we seek for probabilistically the so called pathwise/strong solutions. While in sense of PDE, we consider both the mild/variational solutions and the weak solutions under two different sets of conditions. From now on, the term of weak  solutions are reserved for the weak solutions in the sense of PDE.  The paper is divided into two parts. In the first part, we establish the existence and uniqueness of mild/variational solutions to system (\ref{eq system 00}). To this end, we first appropriately cut off the coefficients of the system and construct a local (in time) mild/variational  solution using fixed point arguments in a certain Banach space and we then show that  the mild/variational solution is global by providing some energy estimates. In the second part, we obtain the existence and uniqueness of pathwise weak solution of the system (\ref{eq system 00}). For this purpose, we first establish the existence of a martingale weak solution. In order to do so, we define a sequence of approximating systems and prove that a subsequence of the approximate solutions converges in law to a martingale weak solution of system (\ref{eq system 00}). Then we  prove that the pathwise uniqueness of weak solutions holds. As an application of Watanable and Yamada Theorem we obtain both the pathwise existence and uniqueness of the weak solution.  Because the proofs of the main results are involved, we will state the main results in next section and leave the details of the arguments in the rest of the paper.
\vskip 0.4cm
The paper is organized as follows. In Section 2, we spread out the precise assumptions and the framework. We also state the main results. Section 3 consists of several subsections. It  is devoted to establishing the existence and uniqueness of mild/variational solution. The entire  Section 4 is to prove the existence and uniqueness of the pathwise weak solution.

\section{Framework and Statement of the Main Results}
\setcounter{equation}{0}
Let $L^q(\mathcal{O})$ denote the $L^q$ space with respect to the Lebesgue measure. $W^{k,q}(\mathcal{O})$ denotes the Sobolev space of functions whose distributional derivatives of  order up to  $k$ belong to $L^q$. Let $A$ be the realization of the Stokes operator $-\mathcal{P}\Delta$, where $\mathcal{P}$ denotes the Helmholtz projection from $L^2(\mathcal{O})$ into the space $H=\{\varphi\in L^2(\mathcal{O})|\nabla\cdot\varphi=0\}$. In the sequel, $\Big(e^{t\Delta}\Big)_{t\geq 0}$, $\Big(e^{-t A}\Big)_{t\geq 0}$ will denote respectively the Neumann heat semigroup and
the Stokes semigroup with Dirichlet boundary condition.

For simplicity, we set $H^k(\mathcal{O}):=W^{k,2}(\mathcal{O})$,
$$
\|\cdot\|_\infty:=\|\cdot\|_{L^\infty(\mathcal{O})},\ \ \ \ \|\cdot\|_{L^q}:=\|\cdot\|_{L^q(\mathcal{O})},\ \ \ \|\cdot\|_\alpha:=\|\cdot\|_{D(A^\alpha)},\ \ \|\cdot\|_{k,q}:=\|\cdot\|_{W^{k,q}(\mathcal{O})},\ \ \|\cdot\|_{H^k}:=\|\cdot\|_{W^{k,2}(\mathcal{O})}.
$$

We introduce the following conditions on the parameters and functions involved in the system (\ref{eq system 00}):
\begin{itemize}
\item[({\bf H.1})] \begin{itemize}
             \item[(a)] $\chi\in C^2([0,\infty)),\ \chi>0$ in $[0,\infty)$,

             \item[(b)] $k\in C^2([0,\infty)),\ k(0)=0,\ k>0$ in $(0,\infty)$,

             \item[(c)] $\phi\in C^2(\bar{\mathcal{O}})$,
           \end{itemize}
\item[({\bf H.2})] $(\frac{k(c)}{\chi(c)})'>0,\ (\frac{k(c)}{\chi(c)})''\leq0,\ (\chi(c)\cdot k(c))'\geq0$ on $[0,\infty)$.
\end{itemize}

Let $U$ be a real Hilbert space and  $\{W_t,\ t\geq0\}$ a $U$-cylindrical  Wiener process
on a given complete, filtered probability space $(\Omega,\mathcal{F},\mathcal{F}_t,\ t\geq0, \mathbb{P})$,
representing the driving  external random force.
\vskip 0.3cm
 Let ${\mathcal{L}^2}(U,D(A^\beta))$ denote the space of Hilbert-Schmidt operators $g$
from $U$ into $D(A^\beta)$ and its norm is denoted by $\|g\|_{{\mathcal{L}^2_\beta}}$. For a  mapping $\sigma:D(A^\beta)\rightarrow \mathcal{L}^2(U,D(A^\beta))$, we introduce
the following hypothesis:
\begin{itemize}
\item[({\bf H.3})]  there exists a positive constant $K$ such that for all $u_1,u_2,u\in H$,
           \begin{eqnarray*}
               \|\sigma(u_1)-\sigma(u_2)\|^2_{\mathcal{L}^2_0}\leq K \|u_1-u_2\|^2_{H},
               \text{ and }
               \|\sigma(u)\|^2_{\mathcal{L}^2_0}\leq K(1+\|u\|^2_{H}),
           \end{eqnarray*}
           where $\mathcal{L}^2_0={\mathcal{L}^2}(U,H)$,
\item[({\bf H.4})]  there exists a positive constant $K$ such that for all $u_1,u_2,u\in D(A^{\alpha})$,
           \begin{eqnarray*}
           \|\sigma(u_1)-\sigma(u_2)\|^2_{\mathcal{L}^2_\alpha}\leq K \|u_1-u_2\|^2_{\alpha},
               \text{ and }
               \|\sigma(u)\|^2_{\mathcal{L}^2_{\alpha}}\leq K(1+\|u\|^2_{{\alpha}}),
           \end{eqnarray*}

\item[({\bf H.5})]  there exists a positive constant $K$ such that for all $u_1,u_2,u\in D(A^{\frac{1}{2}})$,
           \begin{eqnarray*}
           \|\sigma(u_1)-\sigma(u_2)\|^2_{\mathcal{L}^2_{\frac{1}{2}}}\leq K \|u_1-u_2\|^2_{\frac{1}{2}},
               \text{ and }
               \|\sigma(u)\|^2_{\mathcal{L}^2_{{\frac{1}{2}}}}\leq K(1+\|u\|^2_{{\frac{1}{2}}}).
           \end{eqnarray*}

\end{itemize}
 Set:  $u(t)=u(t, \cdot)$, $n(t)=n(t, \cdot)$ and $c(t)=c(t, \cdot)$. Let $q>2$.

\begin{definition}\label{def local solution}
We say that $(n,c,u)$ is a mild solution of system (\ref{eq system 00}) if $(n,c,u)$ is a progressively measurable stochastic processes with values in $C^0(\bar{\mathcal{O}})\times W^{1,q}({\mathcal{O}})\times D(A^\alpha)$, which satisfies
\begin{eqnarray}\label{eq system -1}
&& n(t)=e^{t\delta \Delta}n_0-\int_0^te^{(t-s)\delta \Delta}\Big\{u(s)\cdot \nabla n(s)\Big\}ds-\int_0^te^{(t-s)\delta \Delta}\Big\{\nabla\cdot\Big(\chi(c(s))n(s)\nabla c(s)\Big)\Big\}ds,\nonumber\\
&& c(t)=e^{t\mu \Delta}c_0-\int_0^te^{(t-s)\mu \Delta}\Big\{u(s)\cdot \nabla c(s)\Big\}ds-\int_0^te^{(t-s)\mu \Delta}\Big\{k(c(s))n(s)\Big\}ds,\nonumber\\
&& u(t)=e^{-t\nu A}u_0-\int_0^te^{-(t-s)\nu A}\mathcal{P}\Big\{(u(s)\cdot \nabla)u(s)\Big\}ds-\int_0^te^{-(t-s)\nu A}\mathcal{P}\Big\{n(s)\nabla\phi \Big\}ds \nonumber\\
&&\quad\quad\quad\quad \quad\quad\quad\quad\quad\quad\quad\quad\quad\quad\quad\quad\quad+\int_0^te^{-(t-s)\nu A}\sigma(u(s))dW_s,
\end{eqnarray}
 $P$-a.s.
\end{definition}
\begin{Rem}\label{def local solution-1}
Note that $u(s)\cdot \nabla n(s)=\nabla\cdot(u(s)n(s))$ because $\nabla\cdot u(s)=0$. Under  the setting in the above definition, actually $n(\cdot)\in L^2_{loc}([0, \infty),W^{1,2}({\mathcal{O}}))$ $P$-a.s., and  $(n,c,u)$ is equivalent to a variational solution of the system in the Gelfand triple $W^{1,2}({\mathcal{O}})\subset L^2({\mathcal{O}})\subset W^{1,2}({\mathcal{O}})^*$, that is, $(n,c,u)$  satisfies
\begin{eqnarray}\label{eq system -1}
&& n(t)+\int_0^tu(s)\cdot \nabla n(s)ds=n_0+\delta \int_0^t\Delta n(s)ds-\int_0^t\nabla\cdot\Big(\chi(c(s))n(s)\nabla c(s)\Big)ds,\nonumber\\
&& c(t)+\int_0^tu(s)\cdot \nabla c(s)ds=c_0+\mu \int_0^t\Delta c(s)ds-\int_0^tk(c(s))n(s)ds,\nonumber\\
&& u(t)+\int_0^t\mathcal{P}\Big\{(u(s)\cdot \nabla)u(s)\Big\}ds=u_0-\nu \int_0^tA u(s)ds-\int_0^t\mathcal{P}\Big\{(n(s)\nabla \phi)\Big\}ds\\
&&\quad\quad\quad\quad \quad\quad\quad\quad\quad\quad\quad\quad\quad\quad\quad\quad\quad+\int_0^t\sigma(u(s))dW_s,
\end{eqnarray}
 $P$-a.s..
\end{Rem}
Here is our first main result.
\begin{theorem}\label{Thm main 1}
Assume
\begin{eqnarray}\label{eq boundary condition 3}
&& n_0\in C^0(\bar{\mathcal{O}}),\ \ n_0>0\ \text{in}\ \bar{\mathcal{O}},\nonumber\\
&& c_0\in W^{1,q}(\mathcal{O}),\ \text{for some }q>2,\ c_0>0\ \text{in}\ \bar{\mathcal{O}},\nonumber\\
&& u_0\in D(A^\alpha),\ \text{for some}\ \alpha\in(1/2,1),
\end{eqnarray}
and the assumptions ({\bf H.1})-({\bf H.5}) hold. Then there exists a unique mild/variational solution to the system (\ref{eq system 00}).
\end{theorem}
\vskip 0.4cm

Define $V:=D(A^{1/2})$ and its norm
$$
\|u\|_V:=\|A^{1/2}u\|_H=\|\nabla u\|_{L^2}.
$$
Its dual space will be denoted by $V^*$.

Introduce the following conditions:
\begin{itemize}
\item[({\bf A})] \begin{itemize}
             \item[(a)] $\chi(\cdot)$ and $k(\cdot)$ are smooth with $k(0)=0$, $k(c)>0$ in $(0,\infty)$ and $k'(c)\geq 0$, $\chi(c)>0$ for every $c\in\mathbb{R}$,

             \item[(b)] $\chi'(c)\geq0,\ (\frac{k(c)}{\chi(c)})'>0,\ (\frac{k(c)}{\chi(c)})''<0$, $(\chi(c)\cdot k(c))'>0$ on $[0,\infty)$,

             \item[(c)] $\phi\in C^2(\bar{\mathcal{O}})$,
           \end{itemize}
\item[({\bf B})] $(n_0,c_0,u_0)$ satisfies
         \begin{itemize}
             \item[(B1)] $n_0(x)\geq 0$, $0\leq c_0(x)\leq C_M<\infty$, $\nabla\cdot u_0(x)=0$ on $x\in\mathcal{O}$,

             \item[(B2)] $u_0\in H$,

             \item[(B3)] $n_0(1+|x|+|\ln n_0|)\in L^1(\mathcal{O})$,

             \item[(B4)] $\nabla c_0\in L^2(\mathcal{O})$,
                         $\nabla \Psi (c_0)\in L^2(\mathcal{O})$ where
                         $$
                            \Psi(c)=\int_0^c\sqrt{\frac{\chi(s)}{k(s)}}ds.
                         $$
         \end{itemize}

  \item[({\bf C})] for any $u,u_1,u_2\in H$,
      $$
      \|\sigma(u)\|^2_{\mathcal{L}^2_0}\leq C(1+\|u\|^2_{H})\text{ and }\|\sigma(u_1)-\sigma(u_2)\|^2_{\mathcal{L}^2_0}\leq C\|u_1-u_2\|^2_{H}.
      $$
\end{itemize}
\begin{definition}\label{weak solution}
We say that $({n},{c},{u})$ is a weak  solution to the system (\ref{eq system 00}) if
$({n},{c},{u})$ is a progressively measurable process that satisfies, for any $T>0$,
\begin{itemize}
\item[(1)] $P$-a.s.
\begin{eqnarray*}
 &&{n}(1+|x|+|\ln {n}|)\in L^\infty([0,T],L^1(\mathcal{O})),\ \ \nabla\sqrt{{n}}\in L^2([0,T],L^2(\mathcal{O})),\\
 && {c}\in L^{\infty}([0,T],L^\infty(\mathcal{O})\cap H^1(\mathcal{O}))\cap L^2([0,T],H^2(\mathcal{O})),\\
 && {u}\in C([0,T],H)\cap L^2([0,T],V);
\end{eqnarray*}

\item[(2)]
For all $\psi_1$, $\psi_2\in C^\infty([0,T]\times\mathcal{O})$ with compact supports in the space variable, and $\psi_1(T,\cdot)=\psi_2(T,\cdot)=0$, $P$-a.s.
\begin{eqnarray*}
\int_{\mathcal{O}}\psi_1(0,x)n_0dx
=
\int_0^T\int_{\mathcal{O}}
     {n}[\partial_t\psi_1+\nabla\psi_1\cdot {u}+\delta\Delta\psi_1+\nabla\psi_1\cdot(\chi({c})\nabla {c})]
     dxdt,
\end{eqnarray*}
\begin{eqnarray*}
\int_{\mathcal{O}}\psi_2(0,x)c_0dx
=
\int_0^T\int_{\mathcal{O}}
    {c}[\partial_t\psi_2+\nabla\psi_2\cdot {u}+\mu\Delta\psi_2]-{n}k({c})\psi_2
    dxdt,
\end{eqnarray*}

\item[(3)]
For all $e\in V$, $0\leq t\leq T$,
\begin{eqnarray*}
\langle {u}(t),e\rangle_{H,H}
&=&
\langle u_0,e\rangle_{H,H}-\int_0^t \nu\langle A {u}(s),e\rangle_{V^*,V}ds
-
\int_0^t \langle ({u}(s)\cdot\nabla){u}(s),e\rangle_{V^*,V}ds \\
&&-\int_0^t \langle {n}(s)\nabla \phi,e\rangle_{H,H}ds
 +
 \int_0^t\langle \sigma({u}(s))d{W}_s,e\rangle_{H,H}\\
 \text{holds}&&\quad\quad P-a.s.
\end{eqnarray*}
\end{itemize}
\end{definition}
The following is our second main result.
\begin{theorem}
Assume
the assumptions (A)-(C) hold, and  the function $\chi(\cdot)$ is a positive constant. Then there exists a unique weak solution to the system (\ref{eq system 00}).
\end{theorem}
\vskip 0.4cm

We end this section by recalling the following two  properties of the solution (see Lemma 2.2 in \cite{Winkler}). The first property follows by integrating the first equation in the system (\ref{eq system 00}).
The second one is a consequence of the comparison theorem/maximum principle.
\begin{lemma}\label{lem basic Prop n c}
The solution of (\ref{eq system 00}) satisfies, for all $t\geq 0$,
\begin{eqnarray}\label{eq basic Prop n}
 \int_{\mathcal{O}}n(t,x)dx = \int_{\mathcal{O}}n_0(x)dx,
\end{eqnarray}
and
\begin{eqnarray}\label{eq basic Prop c}
\|c(t,\cdot)\|_{{\infty}}\leq \|c_0\|_{{\infty}},\ \ n(t,x)\geq 0,\ \ c(t,x)\geq 0.
\end{eqnarray}
\end{lemma}
\vskip 0.4cm
Using (\ref{eq basic Prop n}) and the Gagliardo-Nirenberg-Sobolev inequality, we also have
\begin{eqnarray}\label{eq GNS n}
\|n(t)\|_{L^2 } &\leq& C\Big(\|n(t)\|^{1/2}_{L^1 }\|\nabla \sqrt{n(t)}\|_{L^2 }+\|\sqrt{n(t)}\|^2_{L^2 }\Big)\nonumber\\
      &\leq&
         C\Big(\|n_0\|^{1/2}_{L^1 }\|\nabla \sqrt{n(t)}\|_{L^2 }+\|n_0\|_{L^1 }\Big).
\end{eqnarray}

\section{Existence and Uniqueness of Mild/Variational Solutions}
\setcounter{equation}{0}
In this section, we assume that conditions ({\bf H.1})-({\bf H.5}) hold. Our aim is to prove Theorem \ref{Thm main 1}.
\subsection{Existence of Local Solutions}

Introduce the following spaces
$$
\Upsilon^n_t:=L^\infty([0,t],C^0(\bar{\mathcal{O}})),\ \Upsilon^c_t:=L^\infty([0,t],W^{1,q}(\mathcal{O})),
\ \Upsilon^u_t:=L^\infty([0,t],D(A^\alpha))
$$
with the corresponding norms given by
$$
\|n\|_{\Upsilon^n_t}=\sup_{s\in[0,t]}\|n(s)\|_{\infty},\ \|c\|_{\Upsilon^c_t}=\sup_{s\in[0,t]}\|c(s)\|_{1,q},\ \|u\|_{\Upsilon^u_t}=\sup_{s\in[0,t]}\|u(s)\|_{\alpha}.
$$
\begin{definition}\label{def local solution}
We say that $(n,c,u,\tau)$ is a local mild/variational solution of system (\ref{eq system 00}) if

(1) $\tau$ is a stopping time and $(n,c,u)$ is a progressively measurable stochastic processes with values in $C^0(\bar{\mathcal{O}})\times W^{1,q}(\mathcal{O})\times D(A^\alpha)$,

(2) there exists a nondecreasing sequence of stopping times $\{\tau_l,l\geq 1\}$ with $\tau_l\uparrow\tau$
a.s. as $l\uparrow\infty$, such that  $\{(n(t\wedge\tau_l),c(t\wedge\tau_l),u(t\wedge\tau_l)), t\geq 0\}$ is a mild/variational solution to system (\ref{eq system 00}).
\end{definition}

\begin{theorem}
There exists a local mild/variational solution to the system (\ref{eq system 00}).
\end{theorem}
{\bf Proof}. To use a cut off argument, we will modify the coefficients  in system (\ref{eq system 00}).
Fix a function $\theta\in C^2([0,\infty),[0,1])$ such that
\begin{itemize}
\item[(1)] $\theta(r)=1$, $r\in[0,1]$,

\item[(2)] $\theta(r)=0$, $r>2$,

\item[(3)] $\sup_{r\in[0,\infty)}|\theta'(r)|\leq C<\infty$.
\end{itemize}
Set $\theta_m(\cdot)=\theta(\frac{\cdot}{m})$. For every $m\geq 1$, consider the following system of SPDEs
\begin{eqnarray}\label{eq truc n c}
&& dn+\theta_m(\|u\|_{\Upsilon^u_t})\theta_m(\|n\|_{\Upsilon^n_t})u\cdot \nabla ndt
  =
  \delta \Delta ndt-\theta_m(\|n\|_{\Upsilon^n_t})\theta_m(\|c\|_{\Upsilon^c_t})\nabla\cdot(\chi(c)n\nabla c)dt,\nonumber\\
&& dc+\theta_m(\|u\|_{\Upsilon_t})\theta_m(\|c\|_{\Upsilon^c_t})u\cdot \nabla cdt
  =
  \mu \Delta c dt-\theta_m(\|c\|_{\Upsilon^c_t})\theta_m(\|n\|_{\Upsilon^n_t})k(c)ndt,\nonumber\\
&& du+\theta_m(\|u\|_{\Upsilon_t})(u\cdot \nabla)udt+\nabla Pdt
   =
   \nu \Delta udt-\theta_m(\|n\|_{\Upsilon^n_t})n\nabla\phi dt+\sigma(u)dW_t,\\
&& \nabla\cdot u=0,\ \ \ \ \ t>0,\ x\in\mathcal{O}.\nonumber
\end{eqnarray}
To simplify the exposition, we assume $\delta=\mu=\nu=1$, $\chi(c)=1$, and $k(c)=c$. The general case is entirely similar.
\vskip 0.3cm

Let $S_T$ be the space of all $\{\mathcal{F}_t\}_{t\in[0,T]}$-adapted, $C^0(\bar{\mathcal{O}})\times W^{1,q}(\mathcal{O})\times D(A^\alpha)$-valued  stochastic processes $(n(t),c(t),u(t)), t\geq 0$ such that
$$
\|(n,c,u)\|^2_{S_T}:=\mathbb{E}\Big(\|n\|^2_{\Upsilon^n_T}\Big)+\mathbb{E}\Big(\|c\|^2_{\Upsilon^c_T}\Big)
                     +
                     \mathbb{E}\Big(\|u\|^2_{\Upsilon^u_T}\Big)<\infty.
$$
Then $S_T$ equipped with the norm $\|\cdot\|_{S_T}$ is a Banach space.

\vskip 0.3cm

We introduce a mapping $\Phi=(\Phi_1,\Phi_2,\Phi_3)$ on $S_T$ by defining
\begin{eqnarray}\label{eq def Phi1}
  \Phi_1(n,c,u)(t)
&:=&
  e^{t\Delta}n_0
  -
  \int_0^te^{(t-s)\Delta}\Big\{\theta_m(\|n\|_{\Upsilon^n_s})\theta_m(\|c\|_{\Upsilon^c_s})\nabla\cdot(n\nabla c)\nonumber\\
  &&\ \ \ \ \ \ \ \ \ \ \ \ \ \ \ \ \ \ \ \ \ \ \ \ \ \ \ +
  \theta_m(\|u\|_{\Upsilon^u_s})\theta_m(\|n\|_{\Upsilon^n_s})\nabla\cdot (un)\Big\}(s)ds,
\end{eqnarray}
\begin{eqnarray}\label{eq def Phi2}
  \Phi_2(n,c,u)(t)
&:=&
  e^{t\Delta}c_0
  -
  \int_0^te^{(t-s)\Delta}\Big\{\theta_m(\|n\|_{\Upsilon^n_s})\theta_m(\|c\|_{\Upsilon^c_s})nc\nonumber\\
      &&\ \ \ \ \ \ \ \ \ \ \ \ \ \ \ \ \ \ \ \ \ \ +
      \theta_m(\|u\|_{\Upsilon^u_s})\theta_m(\|c\|_{\Upsilon^c_s})u\cdot\nabla c\Big\}(s)ds,
\end{eqnarray}
and
\begin{eqnarray}\label{eq def Phi3}
\Phi_3(n,c,u)(t)&:=&e^{-tA}u_0-\int_0^te^{-(t-s)A}\theta_m(\|u\|_{\Upsilon^u_s})\mathcal{P}\{(u(s)\cdot \nabla)u(s)\}ds\nonumber\\
                   &&
                   -\int_0^te^{-(t-s)A}\theta_m(\|n\|_{\Upsilon^n_s})\mathcal{P}\{n(s)\nabla\phi\} ds
                   +
                  \int_0^te^{-(t-s)A}\sigma(u(s))dW_s.
\end{eqnarray}

\vskip 0.3cm

Let $B$ denote the operator $-\Delta+1$ in $L^q(\mathcal{O})$ ($q>2$) equipped with Neumann boundary condition.
Then, for $\beta\in(\frac{1}{q},\frac{1}{2})$,  we have the continuous imbedding $D(B^\beta)\hookrightarrow C^0(\bar{\mathcal{O}})$.
Using a similar argument as that in \cite{Winkler} ( page 325), we have
\begin{eqnarray}\label{eq Phi1 esta1}
  &&\|\Phi_1(n,c,u)(t)\|_\infty\nonumber\\
&\leq&
  \|e^{t\Delta}n_0\|_\infty + \varrho \int_0^t\|B^\beta e^{-(t-s)(B-1)}\theta_m(\|n\|_{\Upsilon^n_s})\theta_m(\|c\|_{\Upsilon^c_s})\nabla \cdot (n\nabla c)\|_{L^q}ds\nonumber\\
  &&+
  \varrho \int_0^t\|B^\beta e^{-(t-s)(B-1)}\Big(\theta_m(\|u\|_{\Upsilon^u_s})\theta_m(\|n\|_{\Upsilon^n_s})\nabla\cdot (un)\Big)\|_{L^q}ds\nonumber\\
&\leq&
  \|n_0\|_\infty + \varrho \int_0^t (t-s)^{-\beta-\frac{1}{2}}\theta_m(\|n\|_{\Upsilon^n_s})\theta_m(\|c\|_{\Upsilon^c_s})\| n\nabla c\|_{L^q}ds
  \nonumber\\
  &&+
  \varrho \int_0^t (t-s)^{-\beta-\frac{1}{2}}\theta_m(\|u\|_{\Upsilon^u_s})\theta_m(\|n\|_{\Upsilon^n_s})\|un\|_{L^q}ds\nonumber\\
&\leq&
   \|n_0\|_\infty + \varrho m^2 \int_0^t (t-s)^{-\beta-\frac{1}{2}}ds
 \nonumber\\
&\leq&
  \|n_0\|_\infty + \varrho m^2T^{\frac{1}{2}-\beta},\ \ \forall t\in[0,T],
\end{eqnarray}
here we have used the continuous imbedding $D(A^\alpha)\hookrightarrow C^0(\bar{\mathcal{O}})$.

Fix any $\gamma\in (\frac{1}{2},1)$,
\begin{eqnarray}\label{eq Phi2 esta1}
  &&\|\Phi_2(n,c,u)(t)\|_{1,q}\nonumber\\
&\leq&
  \|e^{t\Delta}c_0\|_{1,q} + \varrho \int_0^t\|B^\gamma e^{-(t-s)(B-1)}\theta_m(\|n\|_{\Upsilon^n_s})\theta_m(\|c\|_{\Upsilon^c_s})nc\|_{L^q}ds\nonumber\\
  &&+
  \varrho \int_0^t\|B^\gamma e^{-(t-s)(B-1)}\Big(\theta_m(\|u\|_{\Upsilon^u_s})\theta_m(\|c\|_{\Upsilon^c_s})u\cdot\nabla c\Big)\|_{L^q}ds\nonumber\\
&\leq&
  \varrho\|c_0\|_{1,q} + \varrho \int_0^t(t-s)^{-\gamma}\theta_m(\|n\|_{\Upsilon^n_s})\theta_m(\|c\|_{\Upsilon^c_s})\|nc\|_{L^q}ds\nonumber\\
  &&+
  \varrho \int_0^t(t-s)^{-\gamma}\theta_m(\|u\|_{\Upsilon^u_s})\theta_m(\|c\|_{\Upsilon^c_s})\|(u\cdot\nabla c)\|_{L^q}ds\nonumber\\
&\leq&
  \varrho\|c_0\|_{1,q} + \varrho m^2 \int_0^t(t-s)^{-\gamma}ds\nonumber\\
&\leq&
  \varrho\|c_0\|_{1,q} + \varrho m^2 T^{1-\gamma},\ \ \ \forall t\in[0,T].
\end{eqnarray}

For $\Phi_3$, we have
\begin{eqnarray}\label{eq Phi 3-0}
 && \|A^\alpha\Phi_3(n,c,u)(t)\|_{L^2}\nonumber\\
&\leq&
   \|e^{-tA}A^\alpha u_0\|_{L^2} + \int_0^t\|e^{-(t-s)A}A^\alpha \theta_m(\|u\|_{\Upsilon^u_s})\mathcal{P}\{(u(s)\cdot \nabla)u(s)\}\|_{L^2}ds\nonumber\\
   &&+
   \int_0^t\|e^{-(t-s)A}A^\alpha \theta_m(\|n\|_{\Upsilon^n_s})\mathcal{P}\{n(s)\nabla\phi\}\|_{L^2} ds
                   +
                  \|\int_0^te^{-(t-s)A}A^\alpha\sigma(u(s))dW_s\|_{L^2}\nonumber\\
&\leq&
\|u_0\|_{\alpha}+I_1(t)+I_2(t)+I_3(t).
\end{eqnarray}
Noticing
$$
\|(u\cdot\nabla )u\|_{L^2}\leq\|u\|_{\infty}\|\nabla u\|_{L^2} \leq C \|A^\alpha u\|^2_{L^2},
$$
we have
\begin{eqnarray}\label{eq Phi 3-1}
 I_1(t)
\leq
 C\int_0^t(t-s)^{-\alpha}\theta_m(\|u\|_{\Upsilon^u_s})\|A^\alpha u(s)\|^2_{L^2}ds
 \leq
 Cm^2\int_0^t(t-s)^{-\alpha}ds
 \leq
 Cm^2t^{1-\alpha}.
\end{eqnarray}
For $I_2$, we have
\begin{eqnarray}\label{eq Phi3 I2}
I_2(t)
\leq
\int_0^t(t-s)^{-\alpha}\theta_m(\|n\|_{\Upsilon^n_s})\|n(s)\nabla\phi\|_{L^2}ds
\leq
m\|\nabla\phi\|_{\infty}\int_0^t(t-s)^{-\alpha}ds
\leq
Cmt^{1-\alpha}.
\end{eqnarray}

To estimate $I_3$, let $Z(t):=\int_0^te^{-(t-s)A}A^\alpha\sigma(u(s))dW_s$. Then $Z$ is the solution of the evolution equation
\begin{eqnarray*}
&& dZ(t)=-AZ(t)dt + A^\alpha \sigma(u(t))dW_t,\\
&& Z(0)=0.
\end{eqnarray*}
Applying $\rm It\hat{o}'s$ Formula, and then the BDG inequality, we have
\begin{eqnarray*}
&&\mathbb{E}(\sup_{t\in[0,T]}\|Z(t)\|^2_{L^2})+2E[\int_0^T\|Z(t)\|^2_{1/2}dt]\nonumber\\
&\leq&
  2\mathbb{E}\Big(\sup_{t\in[0,T]}|\int_0^t\langle Z(s),A^\alpha \sigma(u(s))dW_s\rangle_{L^2}|\Big)
  +
  \mathbb{E}\Big(\int_0^T\|A^\alpha \sigma(u(s))\|^2_{\mathcal{L}^2_0}ds\Big)\nonumber\\
&\leq&
 1/2 \mathbb{E}(\sup_{t\in[0,T]}\|Z(t)\|^2_{L^2}) + C\mathbb{E}\Big(\int_0^T\|A^\alpha \sigma(u(s))\|^2_{\mathcal{L}^2_0}ds\Big)\nonumber\\
&\leq&
 1/2 \mathbb{E}(\sup_{t\in[0,T]}\|Z(t)\|^2_{L^2}) + C\mathbb{E}\Big(\int_0^T 1+ \|A^\alpha u(s)\|^2_{L^2}ds\Big)\nonumber\\
&\leq&
 1/2 \mathbb{E}(\sup_{t\in[0,T]}\|Z(t)\|^2_{L^2}) + CT\Big( 1+ \mathbb{E}\Big(\|u\|^2_{\Upsilon^u_T}\Big)\Big),
\end{eqnarray*}
here we have used Assumption (H.4).

Hence
\begin{eqnarray}\label{eq Phi3 I3}
\mathbb{E}(\sup_{t\in[0,T]}\|I_3(t)\|^2_{L^2})\leq CT\Big( 1+ \mathbb{E}(\|u\|^2_{\Upsilon_T^u})\Big).
\end{eqnarray}

Combining (\ref{eq Phi 3-0})--(\ref{eq Phi3 I3}), we get
\begin{eqnarray}\label{eq Phi3 esta1}
\mathbb{E}(\|\Phi_3(n,c,u)\|^2_{\Upsilon_T^u})
\leq
C\Big(\|u_0\|^2_\alpha+m^4 T^{2-2\alpha} + T \mathbb{E}\Big(\|u\|^2_{\Upsilon_T^u}\Big) +T\Big).
\end{eqnarray}

(\ref{eq Phi1 esta1}), (\ref{eq Phi2 esta1}) and (\ref{eq Phi3 esta1}) together show that $\Phi$ maps $S_T$ into itself.
\vskip 0.3cm

Next we will prove that if $T>0$ is small enough then $\Phi$ can be made  a contraction on $S_T$.
\vskip 0.3cm

Let $(n_1,c_1,u_1), (n_2,c_2,u_2)\in S_T$. We have
\begin{eqnarray}\label{eq Phi1 1-2 01}
&&    \|\Phi_1(n_1,c_1£¬u_1)(t)-\Phi_1(n_2,c_2,u_2)(t)\|_\infty\nonumber\\
&\leq&
      \varrho \|\int_0^te^{(t-s)\Delta}\nabla\cdot
            \Big(\theta_m(\|n_1\|_{\Upsilon^n_s})\theta_m(\|c_1\|_{\Upsilon^c_s})n_1\nabla c_1
                 -
                 \theta_m(\|n_2\|_{\Upsilon^n_s})\theta_m(\|c_2\|_{\Upsilon^c_s})n_2\nabla c_2\Big)ds\|_\infty\nonumber\\
      &&+
      \varrho
      \|\int_0^te^{(t-s)\Delta}\nabla\cdot
          \Big(
             \theta_m(\|u_1\|_{\Upsilon^u_s})\theta_m(\|n_1\|_{\Upsilon^n_s})u_1 n_1-\theta_m(\|u_2\|_{\Upsilon^u_s})\theta_m(\|n_2\|_{\Upsilon^n_s})u_2 n_2
          \Big)
      ds\|_\infty\nonumber\\
&:=&
   I_1(t)+I_2(t),
\end{eqnarray}
where $\varrho$ is some generic constant.
Similar to the proof of (\ref{eq Phi1 esta1}), we have
\begin{eqnarray}\label{eq Phi1 1-2 I1}
  &&I_1(t)\\
&\leq&
  \varrho \int_0^t(t-s)^{-\beta-\frac{1}{2}}\|\theta_m(\|n_1\|_{\Upsilon^n_s})\theta_m(\|c_1\|_{\Upsilon^c_s})n_1\nabla c_1
                 -
                 \theta_m(\|n_2\|_{\Upsilon^n_s})\theta_m(\|c_2\|_{\Upsilon^c_s})n_2\nabla c_2\|_{L^q}ds.\nonumber
\end{eqnarray}
Set
$$
J(s)=\|\theta_m(\|n_1\|_{\Upsilon^n_s})\theta_m(\|c_1\|_{\Upsilon^c_s})n_1(s)\nabla c_1(s)
                 -
                 \theta_m(\|n_2\|_{\Upsilon^n_s})\theta_m(\|c_2\|_{\Upsilon^c_s})n_2(s)\nabla c_2(s)\|_{L^q}.
$$
We will distinguish six cases to bound $J$. By the property of $\theta$ and the Minkowski inequality, we have the following estimates.

(J1) Suppose $\|n_1\|_{\Upsilon^n_s}\vee \|c_1\|_{\Upsilon^c_s}\vee\|n_2\|_{\Upsilon^n_s}\vee\|c_2\|_{\Upsilon^c_s}\leq 2m$. We have
\begin{eqnarray*}
J(s)&=&\|\theta_m(\|n_1\|_{\Upsilon^n_s})\theta_m(\|c_1\|_{\Upsilon^c_s})n_1\nabla c_1
                 -
                 \theta_m(\|n_2\|_{\Upsilon^n_s})\theta_m(\|c_2\|_{\Upsilon^c_s})n_2\nabla c_2\|_{L^q}\\
   &\leq&
      \|\theta_m(\|n_1\|_{\Upsilon^n_s})\theta_m(\|c_1\|_{\Upsilon^c_s})n_1\nabla c_1
                 -
                 \theta_m(\|n_2\|_{\Upsilon^n_s})\theta_m(\|c_2\|_{\Upsilon^c_s})n_1\nabla c_1\|_{L^q}\\
      &&+
      \|\theta_m(\|n_2\|_{\Upsilon^n_s})\theta_m(\|c_2\|_{\Upsilon^c_s})n_1\nabla c_1
                 -
                 \theta_m(\|n_2\|_{\Upsilon^n_s})\theta_m(\|c_2\|_{\Upsilon^c_s})n_2\nabla c_2\|_{L^q}\\
   &\leq& \Big|\theta_m(\|n_1\|_{\Upsilon^n_s})\theta_m(\|c_1\|_{\Upsilon^c_s})
   -\theta_m(\|n_2\|_{\Upsilon^n_s})\theta_m(\|c_2\|_{\Upsilon^c_s})\Big|
          \|n_1(s)\nabla c_1(s)\|_{L^q}\\
          &&+
                 \|n_1\nabla c_1(s)-n_2\nabla c_2(s)\|_{L^q}\\
   &\leq& \Big|\theta_m(\|n_1\|_{\Upsilon^n_s})\theta_m(\|c_1\|_{\Upsilon^c_s})
   -\theta_m(\|n_2\|_{\Upsilon^n_s})\theta_m(\|c_2\|_{\Upsilon^c_s})\Big|
          \|n_1(s)\|_\infty\|\nabla c_1(s)\|_{L^q}\\
          &&+
       \|n_1\|_\infty\|c_1(s)-c_2(s)\|_{1,q}+\|n_1(s)-n_2(s)\|_\infty\|c_2(s)\|_{1,q}\\
   &\leq& 4m^2 \Big(
                    \Big|\theta_m(\|n_1\|_{\Upsilon^n_s})\theta_m(\|c_1\|_{\Upsilon^c_s})-\theta_m(\|n_1\|_{\Upsilon^n_s})\theta_m(\|c_2\|_{\Upsilon^c_s})\Big|\\
                    &&\ \ \ \ \ \ \ \ +
                    \Big|\theta_m(\|n_1\|_{\Upsilon^n_s})\theta_m(\|c_2\|_{\Upsilon^c_s})-\theta_m(\|n_2\|_{\Upsilon^n_s})\theta_m(\|c_2\|_{\Upsilon^c_s})\Big|
               \Big)\\
               &&+
       2m(\|c_1(s)-c_2(s)\|_{1,q}+\|n_1(s)-n_2(s)\|_\infty)\\
   &\leq&
      4m^2 \Big(
                    \Big|\theta_m(\|c_1\|_{\Upsilon^c_s})-\theta_m(\|c_2\|_{\Upsilon^c_s})\Big|
                    +
                    \Big|\theta_m(\|n_1\|_{\Upsilon^n_s})-\theta_m(\|n_2\|_{\Upsilon^n_s})\Big|
               \Big)\\
               &&+
       2m(\|c_1(s)-c_2(s)\|_{1,q}+\|n_1(s)-n_2(s)\|_\infty)\\
   &\leq&
      4m^2\frac{C}{m}\Big(\|c_1-c_2\|_{\Upsilon^c_s}+\|n_1-n_2\|_{\Upsilon^n_s}\Big)
      +
      2m(\|c_1(s)-c_2(s)\|_{1,q}+\|n_1(s)-n_2(s)\|_\infty)\\
   &\leq&
      Cm\Big(\|c_1-c_2\|_{\Upsilon^c_s}+\|n_1-n_2\|_{\Upsilon^n_s}\Big).
\end{eqnarray*}

(J2) Suppose $\|n_1\|_{\Upsilon^n_s}\vee \|c_1\|_{\Upsilon^c_s}> 2m$ and $\|n_2\|_{\Upsilon^n_s}\vee\|c_2\|_{\Upsilon^c_s}> 2m$. We have
\begin{eqnarray*}
J(s)=0.
\end{eqnarray*}

(J3) Suppose $\|n_1\|_{\Upsilon^n_s}>2m$ and $\|n_2\|_{\Upsilon^n_s}\vee\|c_2\|_{\Upsilon^c_s}\leq 2m$.
\begin{eqnarray*}
J(s)&=&\|\theta_m(\|n_2\|_{\Upsilon^n_s})\theta_m(\|c_2\|_{\Upsilon^c_s})n_2\nabla c_2\|_{L^q}\\
   &=&
       |\theta_m(\|n_2\|_{\Upsilon^n_s})-\theta_m(\|n_1\|_{\Upsilon^n_s})|\theta_m(\|c_2\|_{\Upsilon^c_s})\|n_2\nabla c_2\|_{L^q}\\
   &\leq&
       Cm\|n_1-n_2\|_{\Upsilon^n_s}.
\end{eqnarray*}

(J4) Suppose $\|c_1\|_{\Upsilon^n_s}>2m$ and $\|n_2\|_{\Upsilon^n_s}\vee\|c_2\|_{\Upsilon^c_s}\leq 2m$.
\begin{eqnarray*}
J(s)&=&\|\theta_m(\|n_2\|_{\Upsilon^n_s})\theta_m(\|c_2\|_{\Upsilon^c_s})n_2\nabla c_2\|_{L^q}\\
   &=&
       |\theta_m(\|c_2\|_{\Upsilon^c_s})-\theta_m(\|c_1\|_{\Upsilon^c_s})|\theta_m(\|n_2\|_{\Upsilon^n_s})\|n_2\nabla c_2\|_{L^q}\\
   &\leq&
       Cm\|c_1-c_2\|_{\Upsilon^c_s}.
\end{eqnarray*}

The proofs of the following two cases are similar as (J3) and (J4).

(J5) If $\|n_2\|_{\Upsilon^n_s}>2m$ and $\|n_1\|_{\Upsilon^n_s}\vee\|c_1\|_{\Upsilon^c_s}\leq 2m$, then
\begin{eqnarray*}
J(s)\leq Cm\|n_1-n_2\|_{\Upsilon^n_s}.
\end{eqnarray*}

(J6) Suppose $\|c_2\|_{\Upsilon^n_s}>2m$ and $\|n_1\|_{\Upsilon^n_s}\vee\|c_1\|_{\Upsilon^c_s}\leq 2m$. Then
\begin{eqnarray*}
J(s)\leq Cm\|c_1-c_2\|_{\Upsilon^c_s}.
\end{eqnarray*}
Putting (J1)--(J6) together, we get
\begin{eqnarray}\label{eq local Phi1 I1}
J(s)\leq Cm(\|c_1-c_2\|_{\Upsilon^c_s}+\|n_1-n_2\|_{\Upsilon^n_s}).
\end{eqnarray}
Substituting (\ref{eq local Phi1 I1}) into (\ref{eq Phi1 1-2 I1}), we get
\begin{eqnarray}\label{eq Phi1 1-2 I1 1}
I_1(t)
&\leq&
C\varrho m(\|c_1-c_2\|_{\Upsilon^c_T}+\|n_1-n_2\|_{\Upsilon^n_T})\int_0^t(t-s)^{-\beta-\frac{1}{2}}ds\nonumber\\
&\leq&
C\varrho m(\|c_1-c_2\|_{\Upsilon^c_T}+\|n_1-n_2\|_{\Upsilon^n_T})T^{\frac{1}{2}-\beta}.
\end{eqnarray}
Using the similar arguments as  in the proof of  (\ref{eq local Phi1 I1}), we can show
\begin{eqnarray*}
&&\|
   \theta_m(\|u_1\|_{\Upsilon^u_s})\theta_m(\|n_1\|_{\Upsilon^n_s})u_1n_1
   -
   \theta_m(\|u_2\|_{\Upsilon^u_s})\theta_m(\|n_2\|_{\Upsilon^n_s})u_2n_2
\|_{L^q}\\
&\leq&
Cm(\|u_1-u_2\|_{\Upsilon^u_s}+\|n_1-n_2\|_{\Upsilon^n_s}).
\end{eqnarray*}
Thus, similar to  (\ref{eq Phi1 1-2 I1 1}), we have
\begin{eqnarray}\label{eq Phi1 1-2 I2 1}
   I_2(t)
&\leq&
   \varrho \int_0^t(t-s)^{-\beta-\frac{1}{2}}
   \|
   \theta_m(\|u_1\|_{\Upsilon^u_s})\theta_m(\|n_1\|_{\Upsilon^n_s})u_1n_1
   -
   \theta_m(\|u_2\|_{\Upsilon^u_s})\theta_m(\|n_2\|_{\Upsilon^n_s})u_2n_2
  \|_{L^q}
         ds\nonumber\\
&\leq&
   C\varrho m(\|u_1-u_2\|_{\Upsilon^u_T}+\|n_1-n_2\|_{\Upsilon^n_T})T^{\frac{1}{2}-\beta}.
\end{eqnarray}
Substitute (\ref{eq Phi1 1-2 I1 1}) and (\ref{eq Phi1 1-2 I2 1}) into (\ref{eq Phi1 1-2 01}) to get
\begin{eqnarray}\label{eq Phi1 1-2}
&&\|\Phi_1(n_1,c_1,u_1)(t)-\Phi_1(n_2,c_2,u_2)(t)\|_\infty\nonumber\\
&\leq&
C\varrho m\Big(\|c_1-c_2\|_{\Upsilon^c_T}+\|u_1-u_2\|_{\Upsilon^u_T}+\|n_1-n_2\|_{\Upsilon^n_T}\Big)T^{\frac{1}{2}-\beta},
\end{eqnarray}
for $t\leq T$.

\vskip 0.3cm

By a similar reasoning, we can show that
\begin{eqnarray}\label{eq Phi 1-2 03}
  &&  \|\Phi_2(n_1,c_1,u_1)(t)-\Phi_2(n_2,c_2,u_2)(t)\|_{1,q}\\
&\leq&
    C\varrho m\Big(\|c_1-c_2\|_{\Upsilon^c_T}+\|u_1-u_2\|_{\Upsilon^u_T}+\|n_1-n_2\|_{\Upsilon^n_T}\Big)T^{1-\gamma},\nonumber
\end{eqnarray}
for $t\leq T$, here $\gamma$ is a number in $(\frac{1}{2}, 1)$.

\vskip 0.3cm
Now we estimate $\|\Phi_3(n_1,c_1,u_1)-\Phi_3(n_2,c_2,u_2)\|_{\Upsilon^u_T}$. We have
\begin{eqnarray}\label{eq u1-u2 01}
&&    \|\Phi_3(n_1,c_1,u_1)(t)-\Phi_3(n_2,c_2,u_2)(t)\|_\alpha\nonumber\\
&\leq&
    \varrho \int_0^t(t-s)^{-\alpha}\|\theta_m(\|u_1\|_{\Upsilon^u_s})(u_1\cdot \nabla)u_1-\theta_m(\|u_2\|_{\Upsilon^u_s})(u_2\cdot \nabla)u_2\|_{L^2}ds\nonumber\\
    &&+
    \varrho \int_0^t(t-s)^{-\alpha}\|\theta_m(\|n_1\|_{\Upsilon^n_s})n_1\nabla\phi-\theta_m(\|n_2\|_{\Upsilon^n_s})n_2\nabla\phi\|_{L^2}ds\nonumber\\
    &&+
    \varrho \Big\|\int_0^te^{-(t-s)A}A^\alpha(\sigma(u_1)-\sigma(u_2))dW_s\Big\|_{L^2}\nonumber\\
&=&
   \Gamma_1(t)+\Gamma_2(t)+\Gamma_3(t).
\end{eqnarray}

Using the similar arguments as in the proof of (\ref{eq Phi1 1-2 I1 1}), it can be shown that
\begin{eqnarray}\label{eq Gamma 2}
       \Gamma_2(t)
   \leq
       C\varrho \|n_1-n_2\|_{\Upsilon^n_T} \cdot T^{1-\alpha},
\end{eqnarray}
for $t\leq T$.
Note that $\Gamma_3(t)=\|Z(t)\|_{L^2}$, where $Z$ satisfies the following SPDE
\begin{eqnarray*}
&&dZ(s)=-AZ(s)ds+A^\alpha(\sigma(u_1)-\sigma(u_2))dW_s,\\
&&Z(0)=0.
\end{eqnarray*}
Using $\rm It\hat{o}$'s Formula and the BDG inequality,  we have
\begin{eqnarray*}
&&  \mathbb{E}\Big(\sup_{t\in[0,T]}\|Z(t)\|^2_{L^2}\Big)+2\mathbb{E}\Big(\int_0^T\|Z(t)\|^2_{\frac{1}{2}}dt\Big)\nonumber\\
&\leq&
    \mathbb{E}\Big(\int_0^T\|A^\alpha(\sigma(u_1(t))-\sigma(u_2(t)))\|^2_{\mathcal{L}_0^2}dt\Big)\nonumber\\
    &&+
    2\mathbb{E}\Big(\sup_{t\in[0,T]}\Big|\int_0^t\Big\langle Z(s),A^\alpha(\sigma(u_1(s))-\sigma(u_2(s)))dW_s\Big\rangle\Big|\Big)\nonumber\\
&\leq&
    \varrho T \mathbb{E}\Big(\|u_1-u_2\|^2_{\Upsilon^u_T}\Big)
    +
    \frac{1}{2}\mathbb{E}\Big(\sup_{t\in[0,T]}\|Z(t)\|^2_{L^2}\Big),
\end{eqnarray*}
here we have used Assumption (H.4).
Hence,
\begin{eqnarray}\label{eq Gamma 3}
   \mathbb{E}\Big(\sup_{t\in[0,T]}\|\Gamma_3(t)\|^2_{L^2}\Big)
\leq
   \varrho T \mathbb{E}\Big(\|u_1-u_2\|^2_{\Upsilon^u_T}\Big).
\end{eqnarray}

To estimate $\Gamma_1$, set
$$
J_1(s)=\|\theta_m(\|u_1\|_{\Upsilon^u_s})(u_1\cdot \nabla)u_1-\theta_m(\|u_2\|_{\Upsilon^u_s})(u_2\cdot \nabla)u_2\|_{L^2}.
$$
We will bound $J_1$ in four different cases. Set $B(u)=(u\cdot \nabla)u$.
\begin{itemize}
\item[(1)] Suppose $\|u_1\|_{\Upsilon^u_s} \vee \|u_2\|_{\Upsilon^u_s}\leq 2m$.  From  the definition of $\theta_m$, we get
           \begin{eqnarray*}\label{eq case 11}
                  J_1(s)
              &\leq&
                  \|B(u_1)-B(u_2)\|_{L^2} + \Big|\theta_m(\|u_1\|_{\Upsilon^u_s})-\theta_m(\|u_2\|_{\Upsilon^u_s})\Big|\|B(u_2)\|_{L^2}\nonumber\\
              &\leq&
                   \varrho \Big(\|u_1\|_\infty \|\nabla(u_1-u_2)\|_{L^2}+\|u_1-u_2\|_\infty \|\nabla u_2\|_{L^2}\Big)\nonumber\\
                   &&+
                   \varrho \frac{C}{m}\|u_1-u_2\|_{\Upsilon^u_s}\|u_2\|^2_\alpha\nonumber\\
              &\leq&
                   \varrho m \|u_1-u_2\|_{\Upsilon^u_s}.
           \end{eqnarray*}

\item[(2)] Suppose  $\|u_1\|_{\Upsilon^u_s}\leq 2m$ and $\|u_2\|_{\Upsilon^u_s}>2m$.  We have
           \begin{eqnarray*}\label{eq case 21}
                 J_1(s)
              &=&
                 \|\theta_m(\|u_1\|_{\Upsilon^u_s})(u_1\cdot \nabla)u_1\|_{L^2}\nonumber\\
              &=&
                  \Big|\theta_m(\|u_1\|_{\Upsilon^u_s})-\theta_m(\|u_2\|_{\Upsilon^u_s})\Big|\|(u_1\cdot \nabla)u_1\|_{L^2}\nonumber\\
              &\leq&
                 \varrho m \|u_1-u_2\|_{\Upsilon^u_s}.
           \end{eqnarray*}

\item[(3)] Suppose  $\|u_1\|_{\Upsilon^u_s}>2m$ and $\|u_2\|_{\Upsilon^u_s}\leq 2m$. Similar to case (2), we have
           \begin{eqnarray*}\label{eq case 31}
                 J_1(s)
              &\leq&
                \varrho m \|u_1-u_2\|_{\Upsilon^u_s}.
           \end{eqnarray*}

\item[(4)] Suppose  $\|u_1\|_{\Upsilon^u_s} \wedge \|u_2\|_{\Upsilon^u_s}> 2m$. Then
           \begin{eqnarray*}\label{eq case 41}
           J_1(s)=0.
           \end{eqnarray*}
\end{itemize}
Hence, it follows that for all the cases,
\begin{eqnarray}\label{eq Gamma 1}
    \Gamma_1(t)
\leq
    \varrho m \|u_1-u_2\|_{\Upsilon_t} t^{1-\alpha}.
\end{eqnarray}
Combining (\ref{eq u1-u2 01}) (\ref{eq Gamma 2}) (\ref{eq Gamma 3}) and (\ref{eq Gamma 1}) together we arrive at
\begin{eqnarray}\label{eq u1-u2}
    &&\mathbb{E}\Big(\|\Phi_3(n_1,c_1,u_1)-\Phi_3(n_2,c_2,u_2)\|^2_{\Upsilon_T}\Big)\nonumber\\
&\leq&
    \varrho T^{2-2\alpha}\mathbb{E}\Big(\|n_1-n_2\|^2_{\Upsilon^n_T}\Big)
    +
    \varrho(T+m^2T^{2-2\alpha})\mathbb{E}\Big(\|u_1-u_2\|^2_{\Upsilon^u_T}\Big).
\end{eqnarray}
\vskip 0.3cm

 By virtue of (\ref{eq Phi1 1-2}) (\ref{eq Phi 1-2 03}) and (\ref{eq u1-u2}), one can find constants $\rho,\ C_m>0$ such that
\begin{eqnarray}\label{eq fix point}
    \|\Phi(n_1,c_1,u_1)-\Phi(n_2,c_2,u_2)\|^2_{S_T}
\leq
    C_{m}T^\rho \|(n_1,c_1,u_1)-(n_2,c_2,u_2)\|^2_{S_T}.
\end{eqnarray}

Choose $T=T_{m}$ such that $C_{m}T^\rho_m=\frac{1}{2}$.  Then $\Phi$ is a contraction on the space $S_{T_m}$. Applying  the Banach fixed point theorem, we conclude that there exists a unique element $(n_m,c_m,u_m)\in S_{T_m}$ such that $(n_m,c_m,u_m)$ is a solution of (\ref{eq truc n c}) for $t\in[0,T_m]$.

\vskip 0.3cm

Let $S^1_T$ be the space of all $\{\mathcal{F}_t\}_{t\in[0,T]}$-adapted, $C^0(\bar{\mathcal{O}})\times W^{1,q}(\mathcal{O})\times D(A^\alpha)$-valued  stochastic processes $(n(t),c(t),u(t)), t\geq 0$ such that
$$
\|(n,c,u)\|^2_{S^1_T}:=\mathbb{E}\Big(\|n\|^2_{\Upsilon^n_T}\Big)+\mathbb{E}\Big(\|c\|^2_{\Upsilon^c_T}\Big)
                     +
                     \mathbb{E}\Big(\|u\|^2_{\Upsilon^u_T}\Big)<\infty,
$$
and
$$
(n,c,u)=(n_m,c_m,u_m)\text{ on }\quad [0,T_m], \quad P\text{-a.s.}.
$$
Then $(S^1_T,\|\cdot\|_{S^1_T})$ is a Banach space.

\vskip 0.3cm
We introduce a mapping $\Phi^1=(\Phi^1_1,\Phi^1_2,\Phi^1_3)$ on $S^1_T$ by defining
\begin{eqnarray*}
  \Phi^1_1(n,c,u)(t+T_m)
&:=&
  e^{t\Delta}n_m(T_m)
  -
  \int_{T_m}^{T_m+t}e^{(T_m+t-s)\Delta}\Big\{\theta_m(\|n\|_{\Upsilon^n_s})\theta_m(\|c\|_{\Upsilon^c_s})\nabla\cdot(n\nabla c)\nonumber\\
  &&\ \ \ \ \ \ \ \ \ \ \ \ \ \ \ \ \ \ \ \ \ \ \ \ \ \ \ +
  \theta_m(\|u\|_{\Upsilon^u_s})\theta_m(\|n\|_{\Upsilon^n_s})\nabla\cdot (un)\Big\}(s)ds,
\end{eqnarray*}
\begin{eqnarray*}
  \Phi^1_2(n,c,u)(t+T_m)
&:=&
  e^{t\Delta}c_m(T_m)
  -
  \int_{T_m}^{T_m+t}e^{(T_m+t-s)\Delta}\Big\{\theta_m(\|n\|_{\Upsilon^n_s})\theta_m(\|c\|_{\Upsilon^c_s})nc\nonumber\\
      &&\ \ \ \ \ \ \ \ \ \ \ \ \ \ \ \ \ \ \ \ \ \ +
      \theta_m(\|u\|_{\Upsilon^u_s})\theta_m(\|c\|_{\Upsilon^c_s})u\cdot\nabla c\Big\}(s)ds,
\end{eqnarray*}
and
\begin{eqnarray*}
\Phi^1_3(n,c,u)(t+T_m)&:=&e^{-tA}u_m(t+T_m)-\int_{T_m}^{T_m+t}e^{-(T_m+t-s)A}\theta_m(\|u\|_{\Upsilon^u_s})\mathcal{P}\{(u(s)\cdot \nabla)u(s)\}ds\nonumber\\
                   &&
                   -\int_{T_m}^{T_m+t}e^{-(T_m+t-s)A}\theta_m(\|n\|_{\Upsilon^n_s})\mathcal{P}\{n(s)\nabla\phi\} ds
                   +
                  \int_{T_m}^{T_m+t}e^{-(T_m+t-s)A}\sigma(u(s))dW_s.
\end{eqnarray*}

\vskip 0.3cm

Observe that the constant $T_m$ does not depend on the initial datum. Repeating the above arguments, we can solve (\ref{eq truc n c}) for $t\in[T_m,2T_m]$, $[2T_m, 3T_m]$,... and we finally  obtain a unique solution $(n_m,c_m,u_m)\in S_{T}$ of (\ref{eq truc n c}) for any $T>0$.

\vskip 0.3cm
Define
\begin{eqnarray}\label{eq stopping time 1}
\tau_m=\inf\{t>0,\ \|n_m\|_{\Upsilon^n_t}\vee\|c_m\|_{\Upsilon^c_t}\vee\|u_m\|_{\Upsilon^u_t}\geq m\}.
\end{eqnarray}
The $\tau_m$ is a stopping time.
When $m>> \|n_0\|_\infty \vee \|c_0\|_{1,q} \vee \|u_0\|_\alpha$, we have
$$
\mathbb{P}(\tau_m>0)=1.
$$
By the definition of $\theta_m$, it is seen that $(n_m(t), c_m(t),\ u_m(t))_{t\in[0,\tau_m)}$ is a local variational solution to the system  (\ref{eq system 00}).\hfill$\Box$

\subsection{Uniqueness of Local Solutions}
\begin{theorem}\label{uniqueness}
Suppose that $(n_1,c_1,u_1,\tau^1)$ and $(n_2,c_2,u_2,\tau^2)$ are two local mild/variational solutions of system (1.1).
Set $\tau=\tau^1 \wedge \tau^2$. Then we have
\begin{eqnarray}\label{eq unique 00}
(n_1,c_1,u_1)=(n_2,c_2,u_2)\ \ \ on\  [0,\tau).
\end{eqnarray}
\end{theorem}
{\bf Proof.}
Define
$$
\tau_R^i=\inf\{t\geq 0,\ \|n_i\|_{\Upsilon^n_t}+\|c_i\|_{\Upsilon^c_t}+\|u_i\|_{\Upsilon^u_t}\geq R\}\wedge\tau^i,\ \ \ i=1,2,
$$
and set  $\tau_R=\tau_R^1\wedge \tau_R^2$.

Notice that $\nabla \cdot (u_i(t)n_i(t))=u_i(t)\cdot \nabla n_i(t)$ because $\nabla\cdot u_i=0$, $i=1,2$.  For all $t\in[0,T\wedge\tau_R)$, we have
\begin{eqnarray}\label{eq unique n1-n2 00}
 && d\|n_1(t)-n_2(t)\|^2_{L^2}+2\|\nabla(n_1(t)-n_2(t))\|^2_{L^2}dt\nonumber\\
&=&
   -2\Big\langle u_1(t)\cdot \nabla n_1(t)-u_2(t)\cdot \nabla n_2(t),n_1(t)-n_2(t)\Big\rangle_{L^2} dt\nonumber\\
   &&-
    2\Big\langle \nabla\cdot(n_1(t)\nabla c_1(t))-\nabla\cdot(n_2(t)\nabla c_2(t)),n_1(t)-n_2(t)\Big\rangle_{L^2} dt\nonumber\\
&\leq&
   \frac{1}{2} \|\nabla(n_1(t)-n_2(t))\|^2_{L^2}dt + C \|u_1(t)n_1(t)-u_2(t)n_2(t)\|^2_{L^2}dt\nonumber\\
   &&+
   C\|n_1(t)\nabla c_1(t)-n_2(t)\nabla c_2(t)\|^2_{L^2}dt\nonumber\\
&\leq&
   \frac{1}{2}\|\nabla(n_1(t)-n_2(t))\|^2_{L^2}dt + C \|u_1(t)-u_2(t)\|^2_{L^2}\|n_1(t)\|^2_\infty dt\nonumber\\
   &&+
   C \|n_1(t)-n_2(t)\|^2_{L^2}\|u_2(t)\|^2_\infty dt
   +
   C \|n_1(t)\|^2_\infty\|\nabla (c_1(t)-c_2(t))\|^2_{L^2}dt\nonumber\\
   &&+
   C \|\nabla c_2(t)\|^2_{L^q}\|n_1(t)-n_2(t)\|^2_{L^{\frac{2q}{q-2}}}dt\nonumber\\
&\leq&
   \frac{1}{2}\|\nabla(n_1(t)-n_2(t))\|^2_{L^2}dt
   +
   C_{R} \|n_1(t)-n_2(t)\|^2_{L^{\frac{2q}{q-2}}}dt\nonumber\\
   &&+
   C_{R}\Big(\|u_1(t)-u_2(t)\|^2_{L^2}+\|n_1(t)-n_2(t)\|^2_{L^2}+\|\nabla (c_1(t)-c_2(t))\|^2_{L^2}\Big)dt\nonumber\\
&\leq&
   \|\nabla(n_1(t)-n_2(t))\|^2_{L^2}dt\nonumber\\
   &&+
   C_{R}\Big(\|u_1(t)-u_2(t)\|^2_{L^2}+\|n_1(t)-n_2(t)\|^2_{L^2}+\|\nabla (c_1(t)-c_2(t))\|^2_{L^2}\Big)dt.
\end{eqnarray}
Here for the last inequality, we have used Ehrling's lemma and the compact embedding $W^{1,2}(\mathcal{O})\hookrightarrow L^{\frac{2q}{q-2}}(\mathcal{O})$.

\vskip 0.3cm
Recall that $W^{1,q}(\mathcal{O})$ is continuously embedded into $C^0(\bar{\mathcal{O}})$. For all $t\in[0,T\wedge\tau_R)$,
we have
\begin{eqnarray}\label{eq unique c1-c2 00}
&& d\|c_1(t)-c_2(t)\|^2_{L^2}+2\|\nabla(c_1(t)-c_2(t))\|^2_{L^2}dt\nonumber\\
&=&
   -2\Big\langle u_1(t)\cdot \nabla c_1(t)-u_2(t)\cdot \nabla c_2(t), c_1(t)-c_2(t)\Big\rangle_{L^2} dt\nonumber\\
   &&-
   2\Big\langle n_1(t)c_1(t)-n_2(t)c_2(t), c_1(t)-c_2(t)\Big\rangle_{L^2} dt\nonumber\\
&\leq&
    \|\nabla(c_1(t)-c_2(t))\|^2_{L^2}dt + C \|u_1(t)c_1(t)-u_2(t)c_2(t)\|^2_{L^2}dt\nonumber\\
   &&+
   \|c_1(t)-c_2(t)\|^2_{L^2}dt + \|n_1(t)c_1(t)-n_2(t)c_2(t)\|^2_{L^2}dt\nonumber\\
&\leq&\|\nabla(c_1(t)-c_2(t))\|^2_{L^2}dt
   +
   C \Big[
                    \|u_1(t)\|^2_\infty \|c_1(t)-c_2(t)\|^2_{L^2}
                    +
                    \|c_2(t)\|^2_\infty \|u_1(t)-u_2(t)\|^2_{L^2}
              \Big]dt\nonumber\\
   &&+
   \|c_1(t)-c_2(t)\|^2_{L^2}dt + \Big[
                                     \|n_1(t)\|^2_\infty \|c_1(t)-c_2(t)\|^2_{L^2}
                                       +
                                     \|c_2(t)\|^2_\infty \|n_1(t)-n_2(t)\|^2_{L^2}
                                 \Big]dt\nonumber\\
&\leq& \|\nabla(c_1(t)-c_2(t))\|^2_{L^2}dt
   +
   C_{R}\Big[ \|c_1(t)-c_2(t)\|^2_{L^2} + \|n_1(t)-n_2(t)\|^2_{L^2} + \|u_1(t)-u_2(t)\|^2_{L^2}\Big]dt.
\end{eqnarray}

By $\rm It\hat{o}$'s formula, for all $t\in[0,T\wedge\tau_R)$,
\begin{eqnarray}\label{eq unique u1-u2 00}
&& d\|u_1(t)-u_2(t)\|^2_{L^2}+2\|\nabla(u_1(t)-u_2(t))\|^2_{L^2}dt\nonumber\\
&=&
  -2\Big\langle \mathcal{P}\Big\{(u_1(t)\cdot\nabla)u_1(t)-(u_2(t)\cdot\nabla)u_2(t)\Big\},u_1(t)-u_2(t)\Big\rangle_{L^2} dt\nonumber\\
  &&-
  2\Big\langle \mathcal{P}\{n_1(t)\nabla\phi-n_2(t)\nabla\phi\},u_1(t)-u_2(t)\Big\rangle_{L^2} dt\nonumber\\
  &&+
  2\Big\langle \sigma(u_1(t))-\sigma(u_2(t)), u_1(t)-u_2(t)\Big\rangle_{L^2} dW_t + \|\sigma(u_1(t))-\sigma(u_2(t))\|^2_{\mathcal{L}_0^2}dt\nonumber\\
&\leq&
  \|\nabla(u_1(t)-u_2(t))\|^2_{L^2}dt \nonumber\\
  &&+ C\|u_2(t)\|^4_{\infty}\|u_1(t)-u_2(t)\|^2_{L^2}dt \nonumber\\
  &&+
  C\|u_1(t)-u_2(t)\|^2_{L^2}dt + C\|n_1(t)-n_2(t)\|^2_{L^2}dt\nonumber\\
  &&+
  2\Big\langle (\sigma(u_1(t))-\sigma(u_2(t))), u_1(t)-u_2(t)\Big\rangle_{L^2}dW_t\nonumber\\
&\leq&
  \|\nabla(u_1(t)-u_2(t))\|^2_{L^2}dt
  +
  C_R\|u_1(t)-u_2(t)\|^2_{L^2}dt + C\|n_1(t)-n_2(t)\|^2_{L^2}dt\nonumber\\
  &&+
  2\Big\langle (\sigma(u_1(t))-\sigma(u_2(t))),u_1(t)-u_2(t)\Big\rangle_{L^2} dW_t.
\end{eqnarray}
For the first inequality of (\ref{eq unique u1-u2 00}), we have used
\begin{eqnarray*}
&&    \Big|\Big\langle \mathcal{P}\Big\{(u_1(t)\cdot\nabla)u_1(t)-(u_2(t)\cdot\nabla)u_2(t)\Big\},u_1(t)-u_2(t)\Big\rangle_{L^2}\Big|\\
&\leq&
      \frac{1}{2}\|\nabla (u_1(t)-u_2(t))\|^2_{L^2} + C\|u_2(t)\|^4_{L^4}\|u_1(t)-u_2(t)\|^2_{L^2}\\
&\leq&
     \frac{1}{2}\|\nabla (u_1(t)-u_2(t))\|^2_{L^2} + C\|u_2(t)\|^4_{\infty}\|u_1(t)-u_2(t)\|^2_{L^2}.
\end{eqnarray*}

Combining (\ref{eq unique n1-n2 00}) (\ref{eq unique c1-c2 00}) and (\ref{eq unique u1-u2 00}), we get
\begin{eqnarray}\label{eq unique 01}
&&  \|n_1(t)-n_2(t)\|^2_{L^2}+C_{R}\|c_1(t)-c_2(t)\|^2_{L^2}+\|u_1(t)-u_2(t)\|^2_{L^2}\nonumber\\
&\leq&
    \widetilde{C}_R\int_0^t\Big(\|n_1(t)-n_2(t)\|^2_{L^2}+\|c_1(t)-c_2(t)\|^2_{L^2}+\|u_1(t)-u_2(t)\|^2_{L^2}\Big)ds\nonumber\\
    &&+
    2\Big|\int_0^t\Big\langle \sigma(u_1(s))-\sigma(u_2(s)),u_1(s)-u_2(s)\Big\rangle_{L^2} dW_s\Big|,\ \ \text{for all }t\in[0,T\wedge\tau_R).
\end{eqnarray}

Let
\begin{eqnarray*}
\Theta(T)&:=&\mathbb{E}\Big(\sup_{t\in[0,T\wedge\tau_R)}\|n_1(t)-n_2(t)\|^2_{L^2}\Big)
     +
     C_{R}\mathbb{E}\Big(\sup_{t\in[0,T\wedge\tau_R)}\|c_1(t)-c_2(t)\|^2_{L^2}\Big)\\
     &&+
     \mathbb{E}\Big(\sup_{t\in[0,T\wedge\tau_R)}\|u_1(t)-u_2(t)\|^2_{L^2}\Big).
\end{eqnarray*}
Apply  BDG inequality, Assumption (H.3) and Gronwall's lemma to conclude from (\ref{eq unique 01}) that
$$
\Theta(T)=0.
$$
We obtain the uniqueness by noting $\tau_R\uparrow\tau$ as $R\uparrow\infty$.  $\Box$

\subsection{Global Existence}
\begin{definition}\label{def maximal solution}
Let $(n,c,u,\tau)$ be a local mild/variational solution of system (\ref{eq system 00}). If $\limsup_{t\nearrow\tau}(\|n\|_{\Upsilon^n_t}+\|c\|_{\Upsilon^c_t}+\|u\|_{\Upsilon^u_t})=\infty$ on $\{\omega,\,\tau<\infty\}$ a.s.,
then the local  solution $(n,c,u,\tau)$ is called a maximal local  solution.

\end{definition}

Recall the stopping times $\{\tau_m,\, m\in\mathbb{N}\}$ defined in (\ref{eq stopping time 1}). By the uniqueness of local  solution
we proved in Section 3.2, we infer that $\tau_m\leq\tau_{m+1}$ a.s. and
$$
(n_{m+1},c_{m+1},u_{m+1})=(n_m,c_m,u_m)\text{ on }[0,\tau_m).
$$
Introduce a stopping time:
$$
\tau =\lim_{m\rightarrow\infty}\tau_m,
$$
and define a stochastic process $(n,c,u)$ on $[0,\tau )$  by
$$
(n,c,u)=(n_m,c_m,u_m)\text{ on }t\in[0,\tau_m).
$$
Since $\|n_m\|_{\Upsilon^n_{\tau_m}}\vee\|c_m\|_{\Upsilon^c_{\tau_m}}\vee\|u_m\|_{\Upsilon^u_{\tau_m}}\geq m$ on $\{\omega,\,\tau<\infty\}$, we have
$$
\limsup_{t\uparrow\tau }(\|n\|_{\Upsilon^n_t}+\|c\|_{\Upsilon^c_t}+\|u\|_{\Upsilon^u_t})
\geq
\lim_{m\uparrow\infty}(\|n_m\|_{\Upsilon^n_{\tau_m}}+\|c_m\|_{\Upsilon^c_{\tau_m}}+\|u_m\|_{\Upsilon^u_{\tau_m}})
=\infty\text{ on }\{\omega,\,\tau<\infty\}.
$$
 Therefore $(n,c,u,\tau)$ is a
maximal local solution of system (\ref{eq system 00}).
\vskip 0.3cm

To obtain the global existence of the solution of the system (\ref{eq system 00}), we will establish some a priori estimates for $(n, c, u)$ in the space
$L^\infty([0,T\wedge\tau); C^0(\bar{\mathcal{O}})\times W^{1,q}(\mathcal{O})\times D(A^{\alpha}))$ for any $T>0$. We first recall the following results from Corollary 4.2 and Lemma 4.5 in \cite{Winkler}.
\begin{lemma}
Let $p>1$ and $r\in [1, \frac{p}{p-1}]$. Then there exists a constant $C_T$ and $C_p$ such that
\begin{equation}\label{3.1}
\int_0^{T\wedge\tau}\|n(t, \cdot)\|_{L^p}^rdt\leq C_T\left (\int_0^{T\wedge\tau}\int_{\mathcal{O}} \frac{|\nabla n(t,x)|^2}{n(t,x)}dxdt +1\right )^{\frac{(p-1)r}{p}},
\end{equation}
and
\begin{equation}\label{3.12}
\int_{\mathcal{O}}n^p(t,x)dx\leq (\int_{\mathcal{O}}n^p(0,x)dx+1)e^{C_p\int_0^t\int_{\mathcal{O}}|\nabla c(s,x)|^4dxds},\ t\in[0,{T\wedge\tau}).
\end{equation}
\end{lemma}
We start with an estimate of the $L^2$ norm of $u$ and $\nabla u$.
\begin{lemma}\label{lem 3.2}
Let $\theta\in (0,1)$. Then there exists a constant $C_{T,\|u_0\|_{L^2}}$ such that
\begin{eqnarray}\label{3.2}
&&E\Big[\sup_{0\leq t< {T\wedge\tau}}\|u(t)\|_{L^2}^4\Big]+E\Big[\Big(\int_0^{T\wedge\tau}\|\nabla u(t)\|_{L^2}^2dt\Big)^2\Big]\nonumber\\
&\leq& C_{T,\|u_0\|_{L^2}}E\Big[\Big (\int_0^{T\wedge\tau}\int_{\mathcal{O}} \frac{|\nabla n(t,x)|^2}{n(t,x)}dxdt +1\Big )^{\theta}\Big].
\end{eqnarray}
Moreover,
\begin{eqnarray}\label{3.7}
E\Big[\int_0^{T\wedge\tau}\int_{\mathcal{O}}|u(t,x)|^4dxdt\Big]
&\leq& C_{T,\|u_0\|_{L^2}}E\Big[\Big (\int_0^{T\wedge\tau}\int_{\mathcal{O}} \frac{|\nabla n(t,x)|^2}{n(t,x)}dxdt +1\Big )^{\theta}\Big].
\end{eqnarray}
\end{lemma}
\vskip 0.3cm
{\bf Proof}. By Ito's formula,
 \begin{eqnarray}\label{3.3}
&& \|u(t)\|^2_{L^2}+2\int_0^t \|\nabla u(s)\|^2_{L^2}ds-\|u_0\|^2_{L^2}\nonumber\\
&=&-2\int_0^t\langle u(s),\mathcal{P}\{ n(s)\nabla\phi\}\rangle_{L^2} ds+2\int_0^t\langle u(s), \sigma(u(s))dW_s\rangle_{L^2} +\int_0^t\|\sigma(u(s))\|^2_{\mathcal{L}^2_0}ds.
\end{eqnarray}
Let $p:=\frac{4}{4-\theta}$. As in the proof of Lemma 4.3 in \cite{Winkler}, writing $p^{\prime}:=\frac{p}{p-1}$, by the Sobolev imbedding $W^{1,2}(\mathcal{O})\hookrightarrow L^{p^{\prime}}(\mathcal{O})$, H\"o{}lder's inequality, we have, for $t\in[0,{T\wedge\tau})$,
\begin{eqnarray}\label{3.4}
&&2\Big|\int_0^t\langle u(s), \mathcal{P}\{n(s)\nabla\phi\}\rangle_{L^2} ds\Big|\nonumber\\
&\leq& c \int_0^t\|u(s)\|_{L^{p'}}\|n(s)\|_{L^{p}}ds\nonumber\\
&\leq& c \int_0^t\|\nabla u(s)\|_{L^{2}}\|n(s)\|_{L^{p}}ds\nonumber\\
&\leq&\int_0^t\|\nabla u(s)\|_{L^{2}}^2ds+c\int_0^t\|n(s)\|_{L^{p}}^2ds\nonumber\\
&\leq& \int_0^t\|\nabla u(s)\|_{L^{2}}^2ds+C_T\left (\int_0^{T\wedge\tau}\int_{\mathcal{O}} \frac{|\nabla n(t,x)|^2}{n(t,x)}dxdt +1\right )^{\frac{\theta}{2}},
\end{eqnarray}
where Lemma 3.1 was used.
Substitute (\ref{3.4}) into (\ref{3.3}) to obtain
\begin{eqnarray}\label{3.5}
&& \sup_{0\leq t< {T\wedge\tau}}\|u(t)\|^2_{L^2}+\int_0^{T\wedge\tau} \|\nabla u(s)\|^2_{L^2}ds\nonumber\\
&\leq&C_{\|u_0\|_{L^2},T}\left (\int_0^{T\wedge\tau}\int_{\mathcal{O}} \frac{|\nabla n(t,x)|^2}{n(t,x)}dxdt +1\right )^{\frac{\theta}{2}} +C\int_0^{T\wedge\tau}(1+\|u(t)\|^2_{L^2})dt\nonumber\\
&&+c\sup_{0\leq t< {T\wedge\tau}}\Big|\int_0^t\langle u(s), \sigma(u(s))dW_s\rangle_{L^2}\Big|.
\end{eqnarray}
Squaring the above inequality and taking expectation, by the BDG inequality and Assumption (H.3), we get
\begin{eqnarray}\label{3.6}
&& E\Big[\sup_{0\leq t< {T\wedge\tau}}\|u(t)\|^4_{L^2}\Big]+E\Big[\Big(\int_0^{T\wedge\tau} \|\nabla u(s)\|^2_{L^2}ds\Big)^2\Big]\nonumber\\
&\leq&C_{\|u_0\|_{L^2},T}E\Big[\left (\int_0^{T\wedge\tau}\int_{\mathcal{O}} \frac{|\nabla n(t,x)|^2}{n(t,x)}dxdt +1\right )^{\theta}\Big] +C_TE\Big[\int_0^{T\wedge\tau}(1+\|u(t)\|^4_{L^2})dt\Big]
\end{eqnarray}
To complete the proof (\ref{3.2}), we  apply the Gronwall's inequality.
\vskip 0.3cm
 By the Gagliardo-Nirenberg inequality we have
\begin{eqnarray}\label{3.8}
&&\int_0^{T\wedge\tau}\int_{\mathcal{O}}|u(t,x)|^4dxdt\nonumber\\
&\leq& C\int_0^{T\wedge\tau}\|\nabla u(t, \cdot)\|_{L^2}^2 \| u(t, \cdot)\|_{L^2}^2 dt\nonumber\\
&\leq& C\sup_{0\leq t< {T\wedge\tau}}\|u(t)\|_{L^2}^4+\Big(\int_0^{T\wedge\tau}\|\nabla u(t)\|_{L^2}^2dt\Big)^2.
\end{eqnarray}
The assertion (\ref{3.7}) follows from (\ref{3.2}).$\Box$

\begin{corollary}\label{cor 3.1}
Let $\theta\in (0,1)$ . The following statements hold:
\begin{equation}\label{3.9}
E\Big[\int_0^{T\wedge\tau}\int_{\mathcal{O}} \frac{|\nabla n(t,x)|^2}{n(t,x)}dxdt +1\Big]<\infty,
\end{equation}
\begin{eqnarray}\label{3.10}
E\Big[\int_0^{T\wedge\tau}\int_{\mathcal{O}}|\nabla c(t,x)|^4dxdt\Big]
&\leq& CE\Big[\Big(\int_0^{T\wedge\tau}\int_{\mathcal{O}} \frac{|\nabla n(t,x)|^2}{n(t,x)}dxdt +1\Big )^{\theta}\Big].
\end{eqnarray}
\end{corollary}
\vskip 0.3cm
\begin{proof} From the proof of Corollary 4.4 in \cite{Winkler}, we know that
\begin{equation}\label{3.11}
\int_0^{T\wedge\tau}\int_{\mathcal{O}} \frac{|\nabla n(t,x)|^2}{n(t,x)}dxdt +\frac{1}{4}\int_0^{T\wedge\tau}\int_{\mathcal{O}}g(c(t,x))|D^2\rho (c(t,x))|^2dxdt\leq C_1+C_2  \int_0^{T\wedge\tau}\int_{\mathcal{O}}|u(t,x)|^4dxdt,
\end{equation}
and
$$\int_0^{T\wedge\tau}\int_{\mathcal{O}}|\nabla c(t,x)|^4dxdt\leq C_3 \int_0^{T\wedge\tau}\int_{\mathcal{O}}g(c(t,x))|D^2\rho (c(t,x))|^2dxdt,$$
where $g(c)=\frac{k(c)}{\chi(c)}$, $\rho (c)=\int_0^c \frac{d\sigma}{g(\sigma)}$.

\vskip 0.3cm
Both (\ref{3.9}) and (\ref{3.10}) now follows from Lemma \ref{lem 3.2}. $\Box$

\end{proof}

\vskip 0.3cm
To proceed, we recall the following inequality obtained in \cite{Winkler}. For any $p>1$,
\begin{equation}\label{3.12}
\int_{\mathcal{O}}n^p(t,x)dx\leq (\int_{\mathcal{O}}n^p(0,x)dx+1)e^{C_p\int_0^t\int_{\mathcal{O}}|\nabla c(s,x)|^4dxds}.
\end{equation}
For $R>0$, define the stopping time $T_R$ by
\begin{eqnarray}\label{3.13}
T_R&=&\inf\{t>0;\ \int_0^t\int_{\mathcal{O}} \frac{|\nabla n(s,x)|^2}{n(s,x)}dxds>R,\ or\quad \int_0^t\int_{\mathcal{O}}|\nabla c(s,x)|^4dxds>R,\nonumber\\
&& \quad \quad \int_0^t\int_{\mathcal{O}} |\nabla u(s,x)|^2dxds>R,\ or \quad \|u(t)\|_{L^2}>R\}.
\end{eqnarray}
Note that $T_R\rightarrow \infty$ a.s. as $R\rightarrow \infty$.
Set $u^R(t,x): =u(t\wedge T_R, x)$, $c^R(t,x): =c(t\wedge T_R, x)$ and $n^R(t,x): =n(t\wedge T_R, x)$. The following result is crucial for establishing the global existence.
\begin{proposition}
For $R>0$ and $T>0$,  there exists some constant $C_{R,T}>0$ such that
\begin{equation}\label{3.14}
E[\sup_{0\leq t\leq T}\|n^R(t, \cdot )\|_{L^{\infty}}]+E[\sup_{0\leq t\leq T}\|\nabla c^R(t, \cdot )\|_{L^{p}}^2]+E[\sup_{0\leq t\leq T}\|A^{\alpha}u^R(t, \cdot )\|_{L^{2}}^2]\leq C_{R,T}.
\end{equation}
\end{proposition}
\vskip 0.3cm
{\bf Proof}. We will prove the proposition along the same lines as in the proof of Theorem 1.1 in \cite{Winkler}. In view of (\ref{3.12}), for any $p>1$ we have
\begin{equation}\label{3.15}
\int_{\mathcal{O}}|n^R|^p(t,x)dx\leq C_{R, p}.
\end{equation}
Applying Ito's formula and following the similar arguments as in the proofs of (4.16) and (4.17) in \cite{Winkler}, we can show that
\begin{eqnarray}\label{3.16}
&&\|\nabla u^R(t)\|_{L^2}^2+\int_0^{t\wedge T_R}\|\Delta u^R(s)\|_{L^2}^2ds\nonumber\\
&\leq& C_R+C_R\int_0^{t\wedge T_R}\|\nabla u^R(s)\|_{L^2}^4ds\nonumber\\
&&+\int_0^{t\wedge T_R}\langle A^{\frac{1}{2}}u^R(s), A^{\frac{1}{2}}\sigma(u^R(s))dW_s\rangle_{L^2,L^2}+\int_0^{t\wedge T_R} \|\sigma(u^R(s))\|_{\mathcal{L}^2_{\frac{1}{2}}}^2ds
\end{eqnarray}
By Gronwall's inequality, it follows that
\begin{eqnarray}\label{3.17}
&&\|\nabla u^R(t)\|_{L^2}^2+\int_0^{t\wedge T_R}\|\Delta u^R(s)\|_{L^2}^2ds\nonumber\\
&\leq& exp\{C_R\int_0^{t\wedge T_R}\|\nabla u^R(s)\|_{L^2}^2ds\}\nonumber\\
&\times& \left [\sup_{0\leq s\leq t}|\int_0^{s\wedge T_R}\langle A^{\frac{1}{2}}u^R(v), A^{\frac{1}{2}}\sigma(u^R(v))dW_v\rangle_{L^2,L^2}|+\int_0^{t\wedge T_R} \|\sigma(u^R(s))\|_{\mathcal{L}^2_{\frac{1}{2}}}^2ds\right ]
\end{eqnarray}
From the definition of $T_R$,
$$exp\{C_R\int_0^{t\wedge T_R}\|\nabla u^R(s)\|_{L^2}^2ds\}\leq e^{C_RR}.$$
Hence, it follows from the Burkholder inequality and (\ref{3.17}) that for any $p>1$,
\begin{eqnarray}\label{3.18}
&&E[\sup_{0\leq s\leq t}\|\nabla u^R(t)\|_{L^2}^{4p}]+E[(\int_0^{t\wedge T_R}\|\Delta u^R(s)\|_{L^2}^2ds)^{2p}]\nonumber\\
&\leq& C_R E[|\int_0^{t}|\langle A^{\frac{1}{2}}u^R(v), A^{\frac{1}{2}}\sigma(u^R(v))\rangle_{L^2,L^2}|^2dv|^p]+C_RE[\int_0^{t} \|\sigma(u^R(s))\|_{\mathcal{L}^2_{\frac{1}{2}}}^{4p}ds]\nonumber\\
&\leq& C_R\left (1+\int_0^tE[\|\nabla u^R(s)\|_{L^2}^{4p}]ds\right ),
\end{eqnarray}
where we have used the fact that $\|A^{\frac{1}{2}}u\|_{L^2}$ is equivalent to $\|\nabla u^R(s)\|_{L^2}$. By Gronwall's inequality, we get from (\ref{3.18}) that
\begin{eqnarray}\label{3.19}
E[\sup_{0\leq s\leq T}\|\nabla u^R(t)\|_{L^2}^{4p}]+E[(\int_0^{T\wedge T_R}\|\Delta u^R(s)\|_{L^2}^2ds)^{2p}]
&\leq& C_{R,p}.
\end{eqnarray}
Next we show that
$$E[\sup_{0\leq t\leq T}\|A^{\alpha}u^R(t, \cdot )\|_{L^{2}}^2]<\infty.$$
By the variation of constants formula we have
\begin{eqnarray}\label{3.20}
 A^{\alpha}u(t)&=& A^{\alpha}e^{-tA}u_0+\int_0^{t}A^{\alpha}e^{-(t-s)A}n(s)\nabla \phi ds+\int_0^{t}A^{\alpha}e^{-(t-s)A}(u(s)\cdot \nabla)u(s)ds\nonumber\\
 &&+ \int_0^{t}A^{\alpha}e^{-(t-s)A}\sigma(u(s))dW_s\nonumber\\
 &:=&A^{\alpha}e^{-tA}u_0+u_1(t)+u_2(t)+u_3(t).
\end{eqnarray}
Clearly,
\begin{equation}\label{3.21}
 \|A^{\alpha}e^{-(t\wedge T_R)A}u_0\|_{L^2}\leq \|A^{\alpha}u_0\|_{L^2}.
\end{equation}
By virtue of (\ref{3.15}), we have
\begin{eqnarray}\label{3.22}
 \|u_1(t\wedge T_R)\|_{L^2}&\leq &\int_0^{t\wedge T_R}\|A^{\alpha}e^{-(t\wedge T_R-s)A}n^R(s)\nabla \phi\|_{L^2} ds\nonumber\\
 &\leq& C\int_0^{t\wedge T_R}(t\wedge T_R-s)^{-\alpha}\|n^R(s)\|_{L^2}ds\leq C(t\wedge T_R)^{1-\alpha}.
\end{eqnarray}
By H\"o{}lder's inequality and Gagliardo-Nirenberg inequality,  it holds that
\begin{eqnarray}\label{3.23}
 \|u_2(t\wedge T_R)\|_{L^2}&\leq& \int_0^{t\wedge T_R}\|A^{\alpha}e^{-(t\wedge T_R-s)A}(u^R(s)\cdot \nabla)u^R(s)\|_{L^2}ds\nonumber\\
 &\leq &C\int_0^{t\wedge T_R}(t\wedge T_R-s)^{-\alpha}\|(u^R(s)\cdot \nabla)u^R(s)\|_{L^2}ds\nonumber\\
 &\leq& C(\int_0^{t\wedge T_R}(t\wedge T_R-s)^{-p^{\prime}\alpha}ds)^{\frac{1}{p'}}(\int_0^{t\wedge T_R}\|(u^R(s)\cdot \nabla)u^R(s)\|_{L^2}^pds)^{\frac{1}{p}}\nonumber\\
 &\leq& Ct^{\frac{1}{p'}-\alpha}(\int_0^{t\wedge T_R}\|\nabla u^R(s)\|_{L^2}^{2p-2}\|\Delta u^R(s)\|_{L^2}^{2}ds)^{\frac{1}{p}}\nonumber\\
 &\leq& C\sup_{0\leq s\leq t}\|\nabla u^R(s)\|_{L^2}^{2-\frac{2}{p}} (\int_0^{t\wedge T_R}\|\Delta u^R(s)\|_{L^2}^{2}ds)^{\frac{1}{p}}.
\end{eqnarray}
This along with (\ref{3.19}) yields
\begin{equation}\label{3.24}
E[\sup_{0\leq t\leq T}\|u_2(t\wedge T_R)\|_{L^2}]\leq C_{T,R},
\end{equation}
where $C_{T,R}$ is some constant.

To estimate $u_3$ in (\ref{3.20}), we notice that $u_3$ satisfies the SPDE:
$$u_3(t)=-\int_0^tAu_3(s)ds+\int_0^tA^{\alpha}\sigma(u(s))dW_s$$
Applying the Ito formula, we get
\begin{eqnarray}\label{3.25}
 &&\|u_3(t\wedge T_R)\|_{L^2}^2+ 2\int_0^{t\wedge T_R}\langle Au_3(s), u_3(s)\rangle_{L^2,L^2} ds\nonumber\\
 &=& 2\int_0^{t\wedge T_R}\langle A^{\alpha}\sigma(u^R(s))dW_s, u_3(s\wedge T_R)\rangle_{L^2,L^2}+\int_0^{t\wedge T_R}\|\sigma(u^R(s))\|_{\mathcal{L}^2_{\alpha}}^2ds
\end{eqnarray}
By Burkholder inequality we get from (\ref{3.25}) that
\begin{eqnarray}\label{3.26}
 &&E[\sup_{0\leq s\leq t}\|u_3(s\wedge T_R)\|_{L^2}^2]\nonumber\\
 &\leq& CE[(\int_0^{t\wedge T_R}\langle A^{\alpha}\sigma(u^R(s)), u_3(s\wedge T_R)\rangle_{L^2,L^2} ^2ds)^{\frac{1}{2}}]\nonumber\\
 &&+ E[\int_0^{t\wedge T_R}\|\sigma(u^R(s))\|_{\mathcal{L}^2_{\alpha}}^2ds]\nonumber\\
 &\leq&\frac{1}{2}E[\sup_{0\leq s\leq t}\|u_3(s\wedge T_R)\|_{L^2}^2]+CE[\int_0^{t\wedge T_R}\|\sigma(u^R(s))\|_{\mathcal{L}^2_{\alpha}}^2ds]\nonumber\\
 &\leq&\frac{1}{2}E[\sup_{0\leq s\leq t}\|u_3(s\wedge T_R)\|_{L^2}^2]+Ct+CE[\int_0^{t}\|A^\alpha u^R(s)\|_{L^2}^2ds],
\end{eqnarray}
which leads to
\begin{equation}\label{3.27}
 E[\sup_{0\leq s\leq t}\|u_3(s\wedge T_R)\|_{L^2}^2]\leq Ct+CE[\int_0^{t}\|A^\alpha u^R(s)\|_{L^2}^2ds].
\end{equation}
Combining (\ref{3.20}), (\ref{3.22}),(\ref{3.24}) and (\ref{3.27}) together we deduce that
\begin{equation}\label{3.28}
 E[\sup_{0\leq s\leq t}\|A^{\alpha}u(s\wedge T_R)\|_{L^2}^2]\leq C+ Ct+CE[\int_0^{t}\|A^{\alpha}u^R(s)\|_{L^2}^2ds].
\end{equation}
An application of Gronwall's inequality yields
\begin{equation}\label{3.29}
E[\sup_{0\leq t\leq T}\|A^{\alpha}u^R(t, \cdot )\|_{L^{2}}^2]<\infty.
\end{equation}
To bound $\|\nabla c^R(t)\|_{L^q}$, we  use the variation of constants formula
\begin{equation}\label{3.30}
c^R(t)=e^{t\wedge T_R\Delta}c_0-\int_0^{t\wedge T_R}e^{(t\wedge T_R-s)\Delta}u^R(s)\cdot \nabla c(s)ds-\int_0^{t\wedge T_R}e^{(t\wedge T_R-s)\Delta}k(c(s))n^R(s)ds
\end{equation}
to obtain
\begin{eqnarray}\label{3.31}
&&\|\nabla c^R(t)\|_{L^q}\nonumber\\
&\leq& C\|\nabla c_0\|_{L^q}+\int_0^{t\wedge T_R}(t\wedge T_R-s)^{-\frac{1}{2}-(\frac{1}{2}-\frac{1}{q})} \|u^R(s)\cdot \nabla c(s)+k(c(s))n^R(s)\|_{L^2}ds\nonumber\\
&\leq& C+\int_0^{t\wedge T_R}(t\wedge T_R-s)^{-\frac{1}{2}-(\frac{1}{2}-\frac{1}{q})} ( \|u^R(s)\|_{L^{\infty}}\|\nabla c^R(s)\|_{L^2}+\|n^R(s)\|_{L^2})ds.
\end{eqnarray}
We note that  $\|n^R(s)\|_{L^2}\leq C_R$ according to (\ref{3.15}). On the other hand, by Lemma 3.4 in \cite{Winkler} and Gagliardo-Nirenberg inequality, we have
\begin{eqnarray}\label{3.32}
\|\nabla c^R(s)\|_{L^2}^2&\leq& c\int_0^{s\wedge T_R}\|u^R(v)\|_{L^4}^4\leq c\int_0^{s\wedge T_R}\|\nabla u^R(v)\|_{L^2}^2\cdot \|u^R(v)\|_{L^2}^2dv\nonumber\\
&\leq& C_R,
\end{eqnarray}
where the last inequality follows from the definition of $T_R$. (\ref{3.31}) and (\ref{3.32}) together yields
\begin{equation}\label{3.33}
\sup_{0\leq t\leq T}\|\nabla c^R(t)\|_{L^q}\leq C_{T, R}
\end{equation}
To estimate $\|n^R(t)\|_{L^{\infty}}$, we fix $\beta\in (\frac{1}{q}, \frac{1}{2})$ and then $r\in (\frac{1}{\beta}, q)$.  As the proof of (4.26) in \cite{Winkler}, we have
\begin{eqnarray}\label{3.34}
\|n^R(t)\|_{L^\infty}&\leq& C\int_0^{t\wedge T_R}(t\wedge T_R-s)^{-\frac{1}{2}-\beta}|n^R\nabla c^R(s)+n^Ru^R(s)\|_{L^r}ds\nonumber\\
&\leq& Ct^{\frac{1}{2}-\beta} \sup_{0\leq s\leq t}\{ \|n^R(s)\|_{L^{\frac{qr}{q-r}}}\|\nabla c^R(s)\|_{L^q}+\|n^R\|_{L^r}\|u^R(s)\|_{L^\infty}\}\nonumber\\
&\leq& C_{R,T}+C_{R,T}\sup_{0\leq s\leq t}\|A^{\alpha}u^R(s)\|_{L^2},
\end{eqnarray}
where (\ref{3.15}),(\ref{3.33}) and the imbedding $D(A^{\alpha})\hookrightarrow L^{\infty}$  have been used. Now, we can conclude from (\ref{3.29}) that

\begin{equation}\label{3.35}
E[\sup_{0\leq t\leq T}\|n^R(t)\|_{L^\infty}]\leq C_{R,T}
\end{equation}
for some constant $C_{R,T}$. The proof is complete.$\Box$

\begin{theorem}\label{thm strong solution in PDE}
Suppose the conditions in Theorem 2.1 are met. Then, the system (1.1) admits a unique global mild/variational solution.
\end{theorem}
\vskip 0.3cm
{\bf Proof}. Let $(n, c, u, \tau )$ be the maximal local solution of system (1.1) obtained in Section 3.3.
From Proposition 3.1 we see that for any $T>0$, $R>0$,
$$\tau \geq T\wedge T_R$$
 Send $R$, $T$ go to infinity to get the global existence. Uniqueness was proved in Section 3.2 $\Box$

\vskip 0.3cm
\begin{Rem}\label{rem strong classical}
We notice that the unique global mild/variational  solution $(n,c,u)$ obtained in Theorem \ref{thm strong solution in PDE} is a weak solution in the sense of  Definition \ref{weak solution}. We only need to verify the statement (1) in Definition \ref{weak solution}.

\vskip 0.2cm
In the proof of Theorem \ref{thm strong solution in PDE}(see (\ref{3.14})), we have, for any $T>0$,
\begin{eqnarray}\label{eq rem 2 01}
\sup_{0\leq t\leq T}\|n(t, \cdot )\|_{{\infty}}+\sup_{0\leq t\leq T}\|\nabla c(t, \cdot )\|_{L^{p}}^2+\sup_{0\leq t\leq T}\|A^{\alpha}u(t, \cdot )\|_{L^{2}}^2<\infty,\ \ \ P\text{-a.s..}
\end{eqnarray}

Theorem \ref{thm strong solution in PDE}, Lemma \ref{lem 3.2} and (\ref{3.9}) imply that
\begin{eqnarray}\label{eq rem 2 02}
E\Big[\sup_{0\leq t\leq {T}}\|u(t)\|_{L^2}^4\Big]+E\Big[\Big(\int_0^{T}\|\nabla u(t)\|_{L^2}^2dt\Big)^2\Big]
<\infty,
\end{eqnarray}
\begin{eqnarray}\label{eq rem 2 03}
E\Big[\int_0^{T}\int_{\mathcal{O}}|u(t,x)|^4dxdt\Big]
<\infty,
\end{eqnarray}
and
\begin{equation}\label{eq rem 2 04}
E\Big[\int_0^T\|\nabla\sqrt{n(t)}\|^2_{L^2}dt\Big]=E\Big[\int_0^{T}\int_{\mathcal{O}} \frac{|\nabla n(t,x)|^2}{n(t,x)}dxdt \Big]<\infty.
\end{equation}
Combining (\ref{eq rem 2 01})--(\ref{eq rem 2 04}) with Lemma \ref{lem basic Prop n c} and the fact that $\mathcal{O}$ is bounded, it is not difficult to deduce that $P$-a.s.
\begin{eqnarray}\label{eq rem 2 05}
 &&{n}(1+|x|)\in L^\infty([0,T],L^1(\mathcal{O})),\ \ \nabla\sqrt{{n}}\in L^2([0,T],L^2(\mathcal{O})),\nonumber\\
 && {c}\in L^{\infty}([0,T],L^\infty(\mathcal{O})\cap H^1(\mathcal{O}))\cap L^2([0,T],H^2(\mathcal{O})),\nonumber\\
 && {u}\in C([0,T],H)\cap L^2([0,T],V).
\end{eqnarray}
For example, to prove ${c}\in L^{\infty}([0,T],L^\infty(\mathcal{O})\cap H^1(\mathcal{O}))\cap L^2([0,T],H^2(\mathcal{O}))$, one can follow
the proof of Lemma \ref{lem 4.2} below.

\vskip 0.3cm
We now estimate $\int_{\mathcal{O}}n|\ln n|dx$. Since
$$
\int_{\mathcal{O}}n \ln ndx = \int_{\mathcal{O}}n|\ln n|dx-2\int_{\mathcal{O}}n \ln{\frac{1}{n}}\mathcal{I}_{n\leq 1}dx,
$$
and, in view of  (\ref{eq basic Prop n}),
\begin{eqnarray*}
0 \leq \int_{\mathcal{O}} n \ln \frac{1}{n}\mathcal{I}_{n\leq 1}dx
&\leq& C\int_{\mathcal{O}}n^{1/2}\mathcal{I}_{n\leq 1}dx\nonumber\\
&\leq&
C_{\mathcal{O}}(\int_{\mathcal{O}}n dx)^{1/2}\nonumber\\
&=&
C_{\mathcal{O}}\|n_0\|^{1/2}_{L^1},
\end{eqnarray*}
it follows that
\begin{eqnarray}\label{eq es n ln n}
0\leq \int_{\mathcal{O}}n |\ln n|dx\leq \int_{\mathcal{O}}n \ln n dx +C.
\end{eqnarray}
Recall $\Psi(c)=\int_0^c\sqrt{\frac{\chi(s)}{k(s)}}ds$ in Assumption (B). Lemma 3.4 in \cite{Winkler} implies that
\begin{eqnarray}\label{eq rem 2 06}
&&\sup_{t\in[0,T]}\Big[\int_{\mathcal{O}}n(t,x) \ln n(t,x) dx+\frac{1}{2}\int_{\mathcal{O}}|\nabla\Psi(c(t,x))|^2dx\Big]\nonumber\\
&\leq&
\int_{\mathcal{O}}n_0(x) \ln n_0(x) dx+\frac{1}{2}\int_{\mathcal{O}}|\nabla\Psi(c_0(x))|^2dx
+
C\int_0^T\int_{\mathcal{O}}|u(t,x)|^4dxdt.
\end{eqnarray}
The assumption on $n_0,c_0$ (see (\ref{eq boundary condition 3})) implies that
\begin{eqnarray}\label{eq rem 2 07}
\int_{\mathcal{O}}n_0(x) \ln n_0(x) dx+\frac{1}{2}\int_{\mathcal{O}}|\nabla\Psi(c_0(x))|^2dx<\infty.
\end{eqnarray}
Putting (\ref{eq rem 2 03}) and (\ref{eq es n ln n})--(\ref{eq rem 2 07}) together,
\begin{eqnarray}\label{eq rem 2 08}
n|\ln n|\in L^\infty([0,T],L^1(\mathcal{O})),\ \ \ P\text{-a.s..}
\end{eqnarray}
\vskip 0.3cm

(\ref{eq rem 2 05}) and (\ref{eq rem 2 08}) show that the statement (1) in Definition \ref{weak solution} holds.

\end{Rem}

\section{Existence and Uniqueness of Weak Solutions}
\setcounter{equation}{0}
In this part, we assume that conditions (A)-(C) introduced in Section 2 are in place. Our aim is to prove the existence and uniqueness of a weak solution to the system (1.1).
\vskip 0.3cm
Because the operator $A$ is positive self-adjoint with compact resolvent, there is a complete orthonormal basis
$\{e_i,\ i\in\mathbb{N}\}$ in $H$ made of eigenvectors of $A$, with corresponding eigenvalues $0<\beta_i\uparrow\infty$, that is
$$
Ae_i=\beta_i e_i,\ \ \ i=1,2,\cdots.
$$
\subsection{Entropy Function}
Let $y\in L^2(\Omega ;\ L^\infty([0,T], H))$ be an adapted process.
Let $(n, c, u)$ be a solution to  the following system
\begin{eqnarray}\label{eq system 01}
&& dn+u\cdot \nabla ndt=\delta \Delta ndt-\nabla\cdot(\chi(c)n\nabla c)dt,\nonumber\\
&& dc+u\cdot \nabla cdt=\mu \Delta c dt-k(c)n dt,\nonumber\\
&& du+(u\cdot \nabla)udt+\nabla Pdt=\nu \Delta udt-n\nabla\phi dt+\sigma(y)dW_t,\\
&& \nabla\cdot u=0,\ \ \ \ \ t>0,\ x\in\mathcal{O}. \nonumber
\end{eqnarray}

 Recall $\Psi(c)=\int_0^c\sqrt{\frac{\chi(s)}{k(s)}}ds$ and $C_M$ is the constant appeared in Condition (B). Set
\begin{eqnarray*}
&&\mathcal{E}(n,c,u)(t)\\
&=&\int_{\mathcal{O}}n(t) \ln n(t)dx +\frac{1}{2}\|\nabla \Psi(c(t))\|^2_{L^2} +2\delta\int_0^t\|\nabla\sqrt{n(s)}\|^2_{L^2}ds+\frac{K}{\nu}\|u(t)\|^2_{L^2}+K\int_0^t\|\nabla u(s)\|^2_{L^2}ds\nonumber\\
&&+
\int_0^t\mu \sum_{i,j}\int_{\mathcal{O}}\Big|\partial_i\partial_j \Psi- \frac{d}{dc}\sqrt{\frac{k(c)}{\chi(c)}}\partial_i\Psi \partial_j\Psi\Big|^2dxds\\
&&+
\int_0^t\lambda_1\mu\int_{\mathcal{O}}\Big|\nabla \Psi\Big|^4dxds + 2\lambda_0\int_0^t\int_{\mathcal{O}}n\Big|\nabla \Psi\Big|^2dxds,
\end{eqnarray*}
here
\begin{eqnarray*}
&&2\lambda_0:=\min_{c\in[0,C_M]}\frac{(\chi(c)k(c))'}{2\chi(c)},\\
&&2\lambda_1:=\min_{c\in[0,C_M]}-\frac{1}{2}\frac{d^2}{dc^2}\Big(\frac{k(c)}{\chi(c)}\Big),
\end{eqnarray*}
are positive by Condition (A).
\vskip 0.3cm

We have the following result.
\begin{proposition}
It holds that
\begin{equation}\label{7.1}
d\mathcal{E}(n,c,u)(t)\leq
Cdt +C\|u(t)\|^2_{L^2}dt + C\langle \sigma(y(t)),u(t)\rangle_{L^2} dW_t + C\|\sigma(y(t))\|^2_{{\mathcal{L}^2_0}}dt.
\end{equation}
\end{proposition}
{\bf Proof.}
Lemma \ref{lem basic Prop n c} and Condition (B) imply that $c$ and $n$ preserve the nonnegativity of the initial data, moreover,
$\|c(t,\cdot)\|_{\infty}\leq C_M$.

Keeping in mind the boundary condition (\ref{eq boundary condition 1}), as (3.5) in \cite{Liu-Lorz},
we can show that
\begin{eqnarray}\label{eq entropy 01}
&&\frac{d}{dt}\int_{\mathcal{O}}n \ln ndx +\frac{1}{2}\frac{d}{dt}\|\nabla \Psi(c)\|^2_{L^2} +4\delta\|\nabla\sqrt{n}\|^2_{L^2}\nonumber\\
&&+
\mu \sum_{i,j}\int_{\mathcal{O}}\Big|\partial_i\partial_j \Psi- \frac{d}{dc}\sqrt{\frac{k(c)}{\chi(c)}}\partial_i\Psi \partial_j\Psi\Big|^2dx
+
\lambda_1\mu\int_{\mathcal{O}}\Big|\nabla \Psi\Big|^4dx + 2\lambda_0\int_{\mathcal{O}}n\Big|\nabla \Psi\Big|^2dx\nonumber\\
&\leq&
K\|\nabla u\|^2_{L^2}.
\end{eqnarray}

\vskip 0.2cm

\vskip 0.2cm

Now applying $\rm It\hat{o}'s$ formula to $\|u\|^2_{L^2}$ and using (\ref{eq GNS n}), it follows that
\begin{eqnarray*}
  &&d\|u\|^2_{L^2} + 2\nu \|\nabla u\|^2_{L^2}dt\nonumber\\
&\leq&
  2\|\nabla \phi\|_{{\infty}} \|n\|_{L^2} \|u\|_{L^2}dt +2\langle \sigma(y)dW_t,u\rangle_{L^2} + \|\sigma(y)\|^2_{{\mathcal{L}^2_0}}dt\nonumber\\
&\leq&
  \delta \frac{\nu}{K}\|\nabla\sqrt{n}\|^2_{L^2}dt + C \|\nabla\phi\|^2_{{\infty}}\|n_0\|_{L^1}\|u\|^2_{L^2}dt
  +
  C\|\nabla\phi\|_{{\infty}}\|n_0\|_{L^1}\|u\|_{L^2}dt\nonumber\\
  &&+
  2\langle \sigma(y)dW_t,u\rangle_{L^2} + \|\sigma(y)\|^2_{{\mathcal{L}^2_0}}dt.
\end{eqnarray*}
Adding the above inequality to (\ref{eq entropy 01}), we obtain
\begin{eqnarray}\label{eq entropy 02}
&&d\int_{\mathcal{O}}n \ln ndx +\frac{1}{2}d\|\nabla \Psi(c)\|^2_{L^2} +2\delta\|\nabla\sqrt{n}\|^2_{L^2}dt+\frac{K}{\nu}d\|u\|^2_{L^2}+K\|\nabla u\|^2_{L^2}dt\nonumber\\
&&+
\mu \sum_{i,j}\int_{\mathcal{O}}\Big|\partial_i\partial_j \Psi- \frac{d}{dc}\sqrt{\frac{k(c)}{\chi(c)}}\partial_i\Psi \partial_j\Psi\Big|^2dxdt
+
\lambda_1\mu\int_{\mathcal{O}}\Big|\nabla \Psi\Big|^4dxdt + 2\lambda_0\int_{\mathcal{O}}n\Big|\nabla \Psi\Big|^2dxdt\nonumber\\
&\leq&
Cdt +C\|u\|^2_{L^2}dt + C\langle \sigma(y),u\rangle_{L^2} dW_t + C\|\sigma(y)\|^2_{{\mathcal{L}^2_0}}dt.
\end{eqnarray}
$\Box$
\subsection{Energy Estimates for Approximating Solutions}\label{subsec 4.2}
In this section, we consider a sequence of approximating solutions and establish
some necessary energy estimates for the proof of the tightness.
\vskip 0.3cm

Let $H_m=\text{span}\{e_1,\cdots,e_m\}$ and define $P_m:H\rightarrow H_m$ as
$$
P_m\,y=\sum_{i=1}^m\langle y,e_i\rangle_{H,H}e_i.
$$
Set $\sigma_m=P_m\sigma$. Then
\begin{eqnarray*}
\|\sigma_m(u)\|^2_{\mathcal{L}^2_{1/2}}
=
\|A^{1/2}\sigma_m(u)\|^2_{\mathcal{L}^2_{0}}
\leq
\beta_m\|\sigma_m(u)\|^2_{\mathcal{L}^2_{0}}
\leq
C\beta_m(1+\|u\|^2_H)
\leq
C\beta_m(1+\frac{1}{\beta_1^2}\|u\|^2_{1/2}).
\end{eqnarray*}
Similarly, we can prove that
\begin{eqnarray*}
&& \|\sigma_m(u_1)-\sigma_m(u_2)\|^2_{\mathcal{L}^2_{1/2}}\leq C_m\|u_1-u_2\|^2_{1/2},\\
&& \|\sigma_m(u)\|^2_{\mathcal{L}^2_\alpha}\leq C_m(1+\|u\|^2_\alpha),\\
&& \|\sigma_m(u_1)-\sigma_m(u_2)\|^2_{\mathcal{L}^2_{\alpha}}\leq C_m\|u_1-u_2\|^2_{\alpha}.
\end{eqnarray*}
This shows that $\sigma_m$ satisfies Conditions \textbf{(H.3)} \textbf{(H.4)} and \textbf{(H.5)}.
\vskip 0.3cm

For any $(n_0,c_0,u_0)$ satisfying Condition $(B)$ in Section 2, it is easy to see that  one can find $(n^m_0,c^m_0,u^m_0)$ such that

(1) $(n^m_0,c^m_0,u^m_0)$ satisfies (\ref{eq boundary condition 3}),

(2)
\begin{eqnarray*}
&& n_0^m(1+|x|+|\ln n_0^m|)\rightarrow n_0(1+|x|+|\ln n_0|)\text{ in }L^1(\mathcal{O}),\\
&& c_0^m\leq C_M\text{ and }c_0^m\rightarrow c_0\text{ in }W^{1,2}(\mathcal{O}),\\
&& \nabla \Psi(c_0^m)\rightarrow \nabla \Psi(c_0)\text{ in }L^2(\mathcal{O}),\\
&& u_0^m\rightarrow u_0\text{ in }L^2(\mathcal{O}).
\end{eqnarray*}

\vskip 0.3cm
By Theorem \ref{thm strong solution in PDE} and Remark \ref{rem strong classical}, we know that there exists an adapted $C^0(\bar{\mathcal{O}})\times W^{1,q}(\mathcal{O})\times D(A^\alpha)$-valued stochastic process $(n^m,c^m,u^m)$ satisfying the following SPDE:
\begin{eqnarray}\label{eq stong system bijin 0}
&& dn^m+u^m\cdot \nabla n^mdt=\delta \Delta n^mdt-\nabla\cdot(\chi(c^m)n^m\nabla c^m)dt,\nonumber\\
&& dc^m+u^m\cdot \nabla c^mdt=\mu \Delta c^m dt-k(c^m)n^m dt,\nonumber\\
&& du^m+(u^m\cdot \nabla)u^mdt+\nabla Pdt=\nu \Delta u^mdt-n^m\nabla\phi dt+\sigma_m(u^m)dW_t,\\
&& \nabla\cdot u^m=0,\ \ \ \ \ t>0,\ x\in\mathcal{O},\nonumber
\end{eqnarray}
with initial value $(n^m_0,c^m_0,u^m_0)$.
\vskip 0.4cm
In the rest of the this section, we will provide a number of estimates for $(n^m,c^m,u^m)$.
\begin{lemma}
There exists a constant $C_T$ independent of $m$ such that
\begin{equation}\label{energy-1}
\mathbb{E}\Big(\sup_{t\in[0,T]}\mathcal{E}(n^m,c^m,u^m)(t)\Big)
\leq
C_T(\int_{\mathcal{O}}n_0 \ln n_0dx +\|\nabla \Psi(c_0)\|^2_{L^2} +\|u_0\|^2_{L^2}
+1).
\end{equation}
\end{lemma}
{\bf Proof.}
By (\ref{eq entropy 02}), we have the following estimates
\begin{eqnarray}\label{eq Strong uniform esta 00}
&&\mathcal{E}(n^m,c^m,u^m)(t)\nonumber\\
&\leq&
\int_{\mathcal{O}}n_0^m \ln n_0^m(t)dx +\frac{1}{2}\|\nabla \Psi(c_0^m)\|^2_{L^2} +\frac{K}{\nu}\|u_0^m\|^2_{L^2}
+
Ct+C\int_0^t\|u^m(s)\|^2_{L^2}ds\nonumber\\
&&+C\int_0^t\langle\sigma_m(u^m(s))dW_s,u^m(s)\rangle_{L^2}
+C\int_0^t\|\sigma_m(u^m(s))\|^2_{\mathcal{L}^2_0}ds.
\end{eqnarray}
In particular, together with (\ref{eq es n ln n}), we have
\begin{eqnarray}\label{eq strong u 0}
&&\frac{K}{\nu}\|u^m(t)\|^2_{L^2}+K\int_0^t\|\nabla u^m(s)\|^2_{L^2}ds\\
&\leq&
\int_{\mathcal{O}}n_0^m \ln n_0^m(t)dx +\frac{1}{2}\|\nabla \Psi(c_0^m)\|^2_{L^2} +\frac{K}{\nu}\|u_0^m\|^2_{L^2}
+
Ct+C\int_0^t\|u^m(s)\|^2_{L^2}ds\nonumber\\
&&+C\int_0^t\langle\sigma_m(u^m(s))dW_s,u^m(s)\rangle_{L^2}
+C\int_0^t\|\sigma_m(u^m(s))\|^2_{\mathcal{L}^2_0}ds.
\end{eqnarray}
By the BDG inequality and the growth condition (C) on $\sigma$,
\begin{eqnarray}\label{eq strong u 1}
&&\mathbb{E}\Big(\sup_{t\in[0,T]}C|\int_0^t\langle\sigma_m(u^m(s))dW_s,u^m(s)\rangle_{L^2}|\Big)\nonumber\\
&\leq&
C\mathbb{E}\Big(\int_0^T\|u^m(s)\|^2_{L^2}\|\sigma_m(u^m(s)\|^2_{\mathcal{L}^2_0})ds\Big)^{1/2}\nonumber\\
&\leq&
1/2\frac{K}{\nu}\mathbb{E}(\sup_{t\in[0,T]}\|u^m(t)\|^2_{L^2})+C\mathbb{E}\int_0^T(1+\|u^m(s)\|^2_{L^2})ds,
\end{eqnarray}
and
\begin{eqnarray}\label{eq strong u 2}
\mathbb{E}\Big(C\int_0^T\|\sigma_m(u^m(s))\|^2_{\mathcal{L}^2_0}ds\Big)
\leq
C\mathbb{E}\int_0^T(1+\|u^m(s)\|^2_{L^2})ds.
\end{eqnarray}
Substituting (\ref{eq strong u 1}) and (\ref{eq strong u 2}) into (\ref{eq strong u 0}), and applying the Gronwall's Lemma, we obtain
\begin{eqnarray}\label{eq strong u}
&&\mathbb{E}\Big(\sup_{t\in[0,T]}\|u^m(t)\|^2_{L^2}\Big)+\mathbb{E}\Big(\int_0^T\|\nabla u^m(s)\|^2_{L^2}ds\Big)\nonumber\\
&\leq&
C(\int_{\mathcal{O}}n_0^m \ln n_0^mdx +\|\nabla \Psi(c_0^m)\|^2_{L^2} +\|u_0^m\|^2_{L^2}
+T)e^{CT}\nonumber\\
&\leq&
C(\int_{\mathcal{O}}n_0 \ln n_0dx +\|\nabla \Psi(c_0)\|^2_{L^2} +\|u_0\|^2_{L^2}
+T)e^{CT}
\end{eqnarray}
(\ref{eq Strong uniform esta 00}) and (\ref{eq strong u}) together imply that
\begin{eqnarray}\label{eq strong uniform esta 01}
\mathbb{E}\Big(\sup_{t\in[0,T]}\mathcal{E}(n^m,c^m,u^m)(t)\Big)
\leq
C_T(\int_{\mathcal{O}}n_0 \ln n_0dx +\|\nabla \Psi(c_0)\|^2_{L^2} +\|u_0\|^2_{L^2}
+1).
\end{eqnarray}
$\Box$
\vskip 0.3cm
\begin{corollary}\label{cor 4.1}
There exists a  constant $C$ independent of $m$ such that

(a) $0\leq n^m(t,x)$ and $c^m(t,x)\in[0,C_M]$, for all $t\geq 0$, $x\in\mathcal{O}$,

(b)
\begin{eqnarray}\label{eq b1}
\mathbb{E}\Big(\sup_{t\in[0,T]}n^m(t)\ln n^m(t)\Big)\leq C,
\end{eqnarray}
and
\begin{eqnarray}\label{eq b2}
\mathbb{E}\Big(\int_0^T\|\nabla \sqrt{n^m(t)}\|^2_{L^2}dt\Big)\leq C.
\end{eqnarray}

(c)
\begin{eqnarray}\label{eq c}
\mathbb{E}\Big(\sup_{t\in[0,T]}n^m(t)|\ln n^m(t)|\Big)\leq C.
\end{eqnarray}

(d)
\begin{eqnarray}\label{eq d}
\mathbb{E}\Big(\int_0^T\|n^m(t)\|^2_{L^2}dt\Big)
\leq
C\mathbb{E}\Big(\int_0^T\|n^m_0\|_{L^1 }\|\nabla \sqrt{n^m}\|^2_{L^2 }+\|n^m_0\|^2_{L^1 }dt\Big)
\leq
C.
\end{eqnarray}

(e)
\begin{eqnarray*}
\mathbb{E}\Big(\sup_{t\in[0,T]}\|u^m(t)\|^2_{L^2}\Big)+\mathbb{E}\Big(\int_0^T\|\nabla u^m(s)\|^2_{L^2}ds\Big)
\leq
C.
\end{eqnarray*}
\end{corollary}
{\bf Proof}.
(a) follows from the comparison theorem, see Lemma \ref{lem basic Prop n c}.
(b) is a consequence of  (\ref{eq strong uniform esta 01}).
(\ref{eq es n ln n}) and (\ref{eq b1}) imply (c).
By (\ref{eq GNS n}) and (\ref{eq b2}), we have
\begin{eqnarray}\label{eq d-1}
\mathbb{E}\Big(\int_0^T\|n^m(t)\|^2_{L^2}dt\Big)
\leq
C\mathbb{E}\Big(\int_0^T\|n^m_0\|_{L^1 }\|\nabla \sqrt{n^m}\|^2_{L^2 }+\|n^m_0\|^2_{L^1 }dt\Big)
\leq
C.
\end{eqnarray}

(e) is the statement  of (\ref{eq strong u}). $\Box$

\vskip 0.3cm
For $N\geq 1$, put
$$
\Omega_N^m:=\Big\{
           \omega:\,
            \sup_{s\in[0,T]}\|u^m(s)\|^2_{L^2} \bigvee \int_0^T\|\nabla u^m(s)\|^2_{L^2}ds \bigvee \int_0^T\|n^m(s)\|^2_{L^2}ds\leq N
            \Big\}.
$$
By Corollary 4.1 and the Chebyshev's inequality, we find that
\begin{eqnarray}\label{7-2}
&&
\mathbb{P}(\Omega_N^m)\nonumber\\
&&\geq
1-\mathbb{P}( \sup_{s\in[0,T]}\|u^m(s)\|^2_{L^2}>N)
 -\mathbb{P}(\int_0^T\|\nabla u^m(s)\|^2_{L^2}ds>N)
 -\mathbb{P}(\int_0^T\|n^m(s)\|^2_{L^2}ds>N)\nonumber\\
&&\geq
1-\frac{3C}{N}.
\end{eqnarray}
\begin{lemma}\label{lem 4.2} We have
\begin{eqnarray}\label{7-3}
\sup_{s\in[0,T]}\|c^m(s)\|^2_{H^1}+\int_0^T\|c^m(s)\|^2_{H^2}ds\leq C_{\mu,T,N,C_M,\|c(0)\|_{H^1}},\ \ \ \omega\in\Omega_N^m,
\end{eqnarray}
where the constant $C_{\mu,T,N,C_M,\|c(0)\|_{H^1}}$ is independent of $m$.
\end{lemma}
{\bf Proof.}

By the chain rule, we have
\begin{eqnarray*}
&&  \|c^m(t)\|^2_{L^2} + 2\mu \int_0^t\|\nabla c^m(s)\|^2_{L^2}ds\\
&=&
    \|c^m(0)\|^2_{L^2} - 2\int_0^t\langle u^m(s)\cdot \nabla c^m(s),c^m(s)\rangle_{L^2}ds
     -
     2\int_0^t\langle k(c^m(s))n^m(s),c^m(s)\rangle_{L^2}ds\\
&\leq&
    \|c(0)\|^2_{L^2} + 2\int_0^t\langle u^m(s)c^m(s),\nabla c^m(s)\rangle_{L^2}ds
    +
    \sup_{r\in[0,C_M]}k^2(r)\int_0^t\|n^m(s)\|^2_{L^2}ds + \int_0^t\|c^m(s)\|^2_{L^2}ds\\
&\leq&
   \|c(0)\|^2_{L^2} + \mu \int_0^t\|\nabla c^m(s)\|^2_{L^2}ds + \frac{C_M^2}{\mu}\int_0^t\|u^m(s)\|^2_{L^2}ds\\
   &&+
   \sup_{r\in[0,C_M]}k^2(r)\int_0^t\|n^m(s)\|^2_{L^2}ds + \int_0^t\|c^m(s)\|^2_{L^2}ds,
\end{eqnarray*}
where we have used (a) of Corollary \ref{cor 4.1} and $u^m(s)\cdot \nabla c^m(s)=\nabla\cdot(u^m(s)c^m(s))$ (due to $\nabla\cdot u^m(s)=0$).
Hence
\begin{eqnarray*}
&&  \|c^m(t)\|^2_{L^2} + \mu \int_0^t\|\nabla c^m(s)\|^2_{L^2}ds\\
&\leq&
    \|c(0)\|^2_{L^2} + \frac{C_M^2}{\mu}\int_0^T\|u^m(s)\|^2_{L^2}ds
    +
   \sup_{r\in[0,C_M]}k^2(r)\int_0^T\|n^m(s)\|^2_{L^2}ds + \int_0^t\|c^m(s)\|^2_{L^2}ds.
\end{eqnarray*}
By the Gronwall's lemma,
\begin{eqnarray}\label{eq estimate cm 01}
&&  \sup_{t\in[0,T]}\|c^m(t)\|^2_{L^2} + \mu \int_0^T\|\nabla c^m(s)\|^2_{L^2}ds\nonumber\\
&\leq&
   e^T
        \cdot
   \Big(
   \|c(0)\|^2_{L^2} + \frac{C_M^2}{\mu}\int_0^T\|u^m(s)\|^2_{L^2}ds
    +
   \sup_{r\in[0,C_M]}k^2(r)\int_0^T\|n^m(s)\|^2_{L^2}ds
   \Big).
\end{eqnarray}

\vskip 0.3cm

Hence we have for $\omega\in \Omega_N^m$,
\begin{eqnarray}\label{eq P1 star estimite cm}
  \sup_{t\in[0,T]}\|c^m(t)\|^2_{L^2} + \mu \int_0^T\|\nabla c^m(s)\|^2_{L^2}ds
\leq
   e^T
        \cdot
   \Big(
   \|c(0)\|^2_{L^2} + \frac{C_M^2}{\mu}TN
    +
   \sup_{r\in[0,C_M]}k^2(r)N
   \Big).
\end{eqnarray}
\vskip 0.3cm

Next we estimate $\|\nabla c^m(t)\|^2_{L^2}$.
Again by the chain rule,
\begin{eqnarray*}
&&  \|\nabla c^m(t)\|^2_{L^2} + 2\int_0^t\langle u^m(s)\cdot\nabla c^m(s),-\Delta c^m(s)\rangle_{L^2}ds\\
&=&
    \|\nabla c^m(0)\|^2_{L^2} - 2\mu\int_0^t\|\Delta c^m(s)\|^2_{L^2}ds
    -
    2\int_0^t\langle k(c^m(s))n^m(s),-\Delta c^m(s)\rangle_{L^2}ds\\
&\leq&
    \|\nabla c(0)\|^2_{L^2} - 2\mu\int_0^t\|\Delta c^m(s)\|^2_{L^2}ds
    +
    \frac{2}{\mu}\sup_{r\in[0,C_M]}k^2(r)\int_0^t\|n^m(s)\|^2_{L^2}ds
    +
    \frac{\mu}{2}\int_0^t\|\Delta c^m(s)\|^2_{L^2}ds.
\end{eqnarray*}
Noticing
$$
\Big|2\int_0^t\langle u^m(s)\cdot\nabla c^m(s),-\Delta c^m(s)\rangle_{L^2}ds\Big|
\leq
\frac{\mu}{2}\int_0^t\|\Delta c^m(s)\|^2_{L^2}ds
+
\frac{2}{\mu}\int_0^t\|u^m(s)\cdot \nabla c^m(s)\|^2_{L^2}ds,
$$
it follows that
\begin{eqnarray}\label{eq star estimate cm}
    &&\|\nabla c^m(t)\|^2_{L^2} + \mu\int_0^t\|\Delta c^m(s)\|^2_{L^2}ds\\
&\leq&
    \|\nabla c(0)\|^2_{L^2}
    +
    \frac{2}{\mu}\int_0^t\|u^m(s)\cdot \nabla c^m(s)\|^2_{L^2}ds
    +
    \frac{2}{\mu}\sup_{r\in[0,C_M]}k^2(r)\int_0^t\|n^m(s)\|^2_{L^2}ds.\nonumber
\end{eqnarray}

\vskip 0.3cm
By the Gagliardo-Nirenberg inequality, we have
\begin{eqnarray}\label{eq P2 01 um}
\|u^m(s)\|^4_{L^4}\leq 2 \|u^m(s)\|^2_{L^2}\|\nabla u^m(s)\|^2_{L^2}.
\end{eqnarray}
Recall also the Gagliardo-Nirenberg-Sobolev inequality:
\begin{eqnarray}\label{eq GNS L4}
\|f\|_{L^4}\leq C(\|\nabla f\|^{1/2}_{L^2}\|f\|^{1/2}_{L^2}+\|f\|_{L^2}).
\end{eqnarray}
Hence we can find a constant $C>0$ such that
\begin{eqnarray}\label{eq P2 02 cm}
    \|\nabla c^m(s)\|_{L^4}
\leq
    C\Big(
        \|c^m(s)\|^{1/2}_{H^2}\|\nabla c^m(s)\|^{1/2}_{L^2} + \|\nabla c^m(s)\|_{L^2}
     \Big),
\end{eqnarray}
here
$$
\|f\|^2_{H^2}=\|f\|^2_{L^2}+\|\nabla f\|^2_{L^2}
                +
             \sum_{i+j=2,\,i,j\geq 0}\int_{\mathcal{O}}\Big|\frac{\partial^i\partial^j f(x_1,x_2)}{\partial x_1^i\partial x_2^j}\Big|^2dx_1dx_2.
$$

According to Proposition 7.2 in \cite{Taylor}(Page 404), for any $f\in H^2$ satisfying the Neumann boundary condition one has
\begin{eqnarray}\label{eq P2 03 estimate H2}
\|f\|^2_{H^2}
\leq
C(\|\Delta f\|^2_{L^2}+ \|f\|^2_{L^2}+\|\nabla f\|^2_{L^2}).
\end{eqnarray}

By (\ref{eq P2 01 um}) (\ref{eq P2 02 cm}) and (\ref{eq P2 03 estimate H2}), for any $\omega\in \Omega_N^m$, we have
\begin{eqnarray}\label{eq P3 star cm 00}
&&    \frac{2}{\mu}\int_0^t\|u^m(s)\cdot \nabla c^m(s)\|^2_{L^2}ds\nonumber\\
&\leq&
      C_\mu\int_0^t \|u^m(s)\|^2_{L^4}\|\nabla c^m(s)\|^2_{L^4}ds\nonumber\\
&\leq&
     C_\mu \int_0^t\|u^m(s)\|_{L^2}\|\nabla u^m(s)\|_{L^2}
              \Big(
                  \|c^m(s)\|_{H^2}\|\nabla c^m(s)\|_{L^2} + \|\nabla c^m(s)\|^2_{L^2}
              \Big)
            ds\nonumber\\
&\leq&
    C_\mu N^{1/2} \int_0^t\|\nabla u^m(s)\|_{L^2}
                                   \Big(\|\Delta c^m(s)\|_{L^2}+\|\nabla c^m(s)\|_{L^2}+\|c^m(s)\|_{L^2}\Big)
                                         \|\nabla c^m(s)\|_{L^2}
                    ds\nonumber\\
                    &&+
      C_\mu N^{1/2} \int_0^t\|\nabla u^m(s)\|_{L^2}\|\nabla c^m(s)\|^2_{L^2}ds\nonumber\\
&\leq&
    \frac{\mu}{2}\int_0^t\|\Delta c^m(s)\|_{L^2}^2ds + C_\mu N\int_0^t\|\nabla u^m(s)\|^2_{L^2}\|\nabla c^m(s)\|^2_{L^2}ds\nonumber\\
    &&+
    C_\mu N^{1/2}\int_0^t\|\nabla u^m(s)\|_{L^2}\|c^m(s)\|_{L^2}\|\nabla c^m(s)\|_{L^2}ds
    +
    C_\mu N^{1/2} \int_0^t\|\nabla u^m(s)\|_{L^2}\|\nabla c^m(s)\|^2_{L^2}ds\nonumber\\
&\leq&
    \frac{\mu}{2}\int_0^t\|\Delta c^m(s)\|_{L^2}^2ds + C_\mu \int_0^t\Big(3N\|\nabla u^m(s)\|^2_{L^2}+1\Big)\|\nabla c^m(s)\|^2_{L^2}ds\nonumber\\
    &&+
    C_\mu t\sup_{s\in[0,t]}\|c^m(s)\|^2_{L^2}.
\end{eqnarray}
For $\omega\in\Omega_N^m$, substituting (\ref{eq P3 star cm 00}) into (\ref{eq star estimate cm}) we obtain
\begin{eqnarray*}
    &&\|\nabla c^m(t)\|^2_{L^2} + \frac{\mu}{2}\int_0^t\|\Delta c^m(s)\|^2_{L^2}ds\\
&\leq&
    \|\nabla c(0)\|^2_{L^2}
    +
    \frac{2}{\mu}\sup_{r\in[0,C_M]}k^2(r)N
    +
    C_\mu \int_0^t\Big(3N\|\nabla u^m(s)\|^2_{L^2}+1\Big)\|\nabla c^m(s)\|^2_{L^2}ds\\
&&    +
     C_\mu T\sup_{s\in[0,T]}\|c^m(s)\|^2_{L^2}.
\end{eqnarray*}
Hence by Gronwall's lemma and (\ref{eq P1 star estimite cm}), it follows that
\begin{eqnarray}\label{eq 4.30}
    &&\|\nabla c^m(t)\|^2_{L^2} + \frac{\mu}{2}\int_0^t\|\Delta c^m(s)\|^2_{L^2}ds\nonumber\\
&\leq&
     e^{C_\mu\int_0^T(3N\|\nabla u^m(s)\|^2_{L^2}+1)ds}
     \cdot
     \Big[
        \|\nabla c(0)\|^2_{L^2}
         +
         \frac{2}{\mu}\sup_{r\in[0,C_M]}k^2(r)N\nonumber\\
         &&\ \ \ \ \ \ \ \ \ \ \ \ \ \ \ \ \ \ \ \ \ \ \ \ \ \ \ \ \ \ \ \  \ +
         C_\mu T e^T
                   \cdot
                       \Big(
                           \|c(0)\|^2_{L^2} + \frac{C_M^2}{\mu}TN
                             +
                           \sup_{r\in[0,C_M]}k^2(r)N
                        \Big)
     \Big]\nonumber\\
&\leq&
     e^{3C_\mu N^2+C_\mu T}\Big[
        \|\nabla c(0)\|^2_{L^2}
         +
         \frac{2}{\mu}\sup_{r\in[0,C_M]}k^2(r)N\nonumber\\
         &&\ \ \ \ \ \ \ \ \ \ \ \ \ \ \ \ \ +
         C_\mu T e^T
                   \cdot
                       \Big(
                           \|c(0)\|^2_{L^2} + \frac{C_M^2}{\mu}TN
                             +
                           \sup_{r\in[0,C_M]}k^2(r)N
                        \Big)
     \Big].
\end{eqnarray}
Thus, in view of (\ref{eq P2 03 estimate H2}) (\ref{eq P1 star estimite cm}) and (\ref{eq 4.30}), we can conclude that
\begin{eqnarray}\label{eq estimate cm H1 H2}
 \sup_{s\in[0,T]}\|c^m(s)\|^2_{H^1}+\int_0^T\|c^m(s)\|^2_{H^2}ds\leq C_{\mu,T,N,C_M,\|c(0)\|_{H^1}},\ \ \ \omega\in\Omega_N^m.
\end{eqnarray}
The constant $C_{\mu,T,N,C_M,\|c(0)\|_{H^1}}$ is independent of $m$.$\Box$
\begin{corollary}
There exists a constant $C_{\mu,T,N,C_M,\|c(0)\|_{H^1}}$ such that for all $\omega\in\Omega_N^m$,

\begin{eqnarray}\label{7-5}
      \int_0^T\|\frac{d c^m(t)}{dt}\|^2_{L^2}dt
&\leq&C_{\mu,T,N,C_M,\|c(0)\|_{H^1}}.
\end{eqnarray}
\end{corollary}

\vskip 0.3cm
\noindent {\bf Proof}. Combining  (\ref{eq P3 star cm 00}) and (\ref{eq estimate cm H1 H2}), for $\omega\in\Omega_N^m$, we have
\begin{eqnarray*}
    \int_0^T\|u^m(s)\cdot\nabla c^m(s)\|^2_{L^2}ds
\leq
    C_{\mu,T,N,C_M,\|c(0)\|_{H^1}}.
\end{eqnarray*}
Hence, for $\omega\in\Omega_N^m$, by (\ref{eq estimate cm H1 H2}), we have
\begin{eqnarray}\label{eq estimate cm dt}
      &&\int_0^T\|\frac{d c^m(t)}{dt}\|^2_{L^2}dt\nonumber\\
&\leq&
      C\Big[
           \int_0^T\|u^m(s)\cdot\nabla c^m(s)\|^2_{L^2}ds
           +
           \mu^2\int_0^T\|\Delta c^m(s)\|^2_{L^2}ds
           +
           \sup_{r\in[0,C_M]}k^2(r)\int_0^T\|n^m(s)\|^2_{L^2}ds
      \Big]\nonumber\\
&\leq&
      C_{\mu,T,N,C_M,\|c(0)\|_{H^1}},\ \ \ \omega\in\Omega_N^m.
\end{eqnarray}
$\Box$
\vskip 0.3cm
\begin{lemma}
There exists a constant $C(T,N)$ such that for all $\omega\in\Omega_N^m$ and $m\geq 1$,
\begin{equation}\label{7-7}
 \int_0^T\|\frac{d n^m(t)}{dt}\|^2_{H^{-3}}dt\leq
         C(T,N),
\end{equation}
and
\begin{equation}\label{7-8}
 \sup_{t\in[0,T]}\|n^m(t)\|^2_{L^2}+\int_0^T\|\nabla n^m(s)\|^2_{L^2}ds
\leq C(T,N).
\end{equation}
\end{lemma}
{\bf Proof.}

We first prove (\ref{7-7}). According to the Sobolev inequalities, we have
$$
\|\nabla \varphi\|_{\infty}\leq C\|\varphi\|_{H^3},\ \ \varphi\in H^3(\mathcal{O}).
$$
Hence, for $\varphi\in H^3(\mathcal{O})$,
\begin{eqnarray*}
     &&|\langle \varphi,\frac{d n^m(t)}{dt}\rangle_{L^2}|\\
&\leq&
     \delta|\langle \Delta \varphi, n^m(t)\rangle_{L^2}|
     +
     |\langle \nabla \varphi, n^m(t)u^m(t)\rangle_{L^2}|
     +
     |\langle \nabla \varphi, \chi(c^m(t))n^m(t)\nabla c^m(t)\rangle_{L^2}|\\
&\leq&
     \delta \|\Delta \varphi\|_{L^2}\|n^m(t)\|_{L^2} + C\|\varphi\|_{H^3}\|n^m(t)\|_{L^2}\|u^m(t)\|_{L^2}
     +
     C\|\varphi\|_{H^3}\|n^m(t)\|_{L^2}\|\nabla c^m(t)\|_{L^2}.
\end{eqnarray*}
Therefore, for $\omega\in \Omega_N^m$,
\begin{eqnarray}\label{eq estimate nm dt}
&&    \int_0^T\|\frac{d n^m(t)}{dt}\|^2_{H^{-3}}dt\nonumber\\
&\leq&
      C\Big[
         \delta\int_0^T\|n^m(t)\|^2_{L^2}dt
         +
         \int_0^T\|n^m(t)\|^2_{L^2}\|u^m(t)\|^2_{L^2}dt
         +
         \int_0^T\|n^m(t)\|^2_{L^2}\|\nabla c^m(t)\|^2_{L^2}dt
      \Big]\nonumber\\
&\leq&
    C\Big[
         \delta N
         +
         N^2
         +
         N\sup_{t\in[0,T]}\|\nabla c^m(t)\|^2_{L^2}
      \Big]\nonumber\\
&\leq&
     C\Big[
         \delta N
         +
         N^2
         +
         C_{\mu,T,N,C_M,\|c(0)\|_{H^1}}
      \Big].
\end{eqnarray}
(\ref{eq estimate cm H1 H2}) has been used in the last inequality. This proves (\ref{7-7}).

\vskip 0.3cm
By the chain rule,
\begin{eqnarray}\label{eq estimate nm L2 01}
&&\|n^m(t)\|^2_{L^2} + 2\delta\int_0^t\|\nabla n^m(s)\|^2_{L^2}ds\nonumber\\
&=&
   \|n^m(0)\|^2_{L^2}+2\int_0^t\langle u^m(s)n^m(s),\nabla n^m(s)\rangle_{L^2}ds
    +
    2\int_0^t\langle \chi(c^m(s))n^m(s)\nabla c^m(t),\nabla n^m(s)\rangle_{L^2}ds\nonumber\\
&\leq&
  \|n_0\|^2_{L^2}+\delta\int_0^t\|\nabla n^m(s)\|^2_{L^2}ds+\frac{2}{\delta}\int_0^t\|u^m(s)n^m(s)\|^2_{L^2}ds\nonumber\\
  &&+
  \frac{2}{\delta}\sup_{r\in[0,C_M]}\chi^2(r)\int_0^t\|n^m(s)\nabla c^m(s)\|^2_{L^2}ds.
\end{eqnarray}
By (\ref{eq P2 01 um}) and (\ref{eq GNS L4}),
\begin{eqnarray}\label{eq estimate nm L2 02}
&&\frac{2}{\delta}\int_0^t\|u^m(s)n^m(s)\|^2_{L^2}ds\nonumber\\
&\leq&
 C_\delta\int_0^t\|u^m(s)\|^2_{L^4}\|n^m(s)\|^2_{L^4}ds\nonumber\\
&\leq&
C_\delta\int_0^t\|u^m(s)\|_{L^2}\|\nabla u^m(s)\|_{L^2}\Big(\|\nabla n^m(s)\|_{L^2}\|n^m(s)\|_{L^2}+\|n^m(s)\|^2_{L^2}\Big)ds\\
&\leq&
\frac{\delta}{4}\int_0^t\|\nabla n^m(s)\|^2_{L^2}ds
+
C_\delta\int_0^t\|n^m(s)\|^2_{L^2}\Big(\|u^m(s)\|^2_{L^2}\|\nabla u^m(s)\|^2_{L^2}+\|u^m(s)\|_{L^2}\|\nabla u^m(s)\|_{L^2}\Big)ds.\nonumber
\end{eqnarray}
In view of (\ref{eq GNS L4}) and (\ref{eq P2 02 cm}), we have
\begin{eqnarray}\label{eq estimate nm L2 03}
&&\frac{2}{\delta}\sup_{r\in[0,C_M]}\chi^2(r)\int_0^t\|n^m(s)\nabla c^m(s)\|^2_{L^2}ds\\
&\leq&
C_\delta \int_0^t\|n^m(s)\|^2_{L^4}\|\nabla c^m(s)\|^2_{L^4}ds\nonumber\\
&\leq&
C_\delta\int_0^t\Big(\|\nabla n^m(s)\|_{L^2}\|n^m(s)\|_{L^2}+\|n^m(s)\|^2_{L^2}\Big)
               \cdot
                \Big(\|c^m(s)\|_{H^2}\|\nabla c^m(s)\|_{L^2}+\|\nabla c^m(s)\|^2_{L^2}\Big)
                ds\nonumber\\
&\leq&
  \frac{\delta}{4}\int_0^t\|\nabla n^m(s)\|^2_{L^2}ds
+
  C_\delta\int_0^t\|n^m(s)\|^2_{L^2}
                       \Big(
                           \|c^m(s)\|^2_{H^2}\|\nabla c^m(s)\|^2_{L^2}+\|\nabla c^m(s)\|^4_{L^2}\nonumber\\
  &&\ \ \ \ \ \ \ \ \ \ \ \ \ \ \ \ \ \ \ \ \ \ \ \ \ \ \ \ \ \ \ \ \ \ \ \ \ \ \ \ \ \ \ \ \ \ \ \ \ \ \ \ +
                           \|c^m(s)\|_{H^2}\|\nabla c^m(s)\|_{L^2}+\|\nabla c^m(s)\|^2_{L^2}
                      \Big)
          ds.\nonumber
\end{eqnarray}
Combining (\ref{eq estimate nm L2 01})--(\ref{eq estimate nm L2 03}), and applying the Gronwall's lemma,
we obtain
\begin{eqnarray*}
\sup_{t\in[0,T]}\|n^m(t)\|^2_{L^2}+\frac{\delta}{4}\int_0^T\|\nabla n^m(s)\|^2_{L^2}ds
\leq
\|n_0\|^2_{L^2}e^{C_\delta\Xi(T)}
\end{eqnarray*}
where
\begin{eqnarray*}
\Xi(T):&=&\int_0^T\Big(\|u^m(s)\|^2_{L^2}\|\nabla u^m(s)\|^2_{L^2}+\|u^m(s)\|_{L^2}\|\nabla u^m(s)\|_{L^2}\\
                                   &&\ \ \ \ \ \ \ \ \ \  +
                                   \|c^m(s)\|^2_{H^2}\|\nabla c^m(s)\|^2_{L^2}+\|\nabla c^m(s)\|^4_{L^2}\\
                                   &&\ \ \ \ \ \ \ \ \ \ +
                                   \|c^m(s)\|_{H^2}\|\nabla c^m(s)\|_{L^2}+\|\nabla c^m(s)\|^2_{L^2}
                                   \Big)ds.
\end{eqnarray*}
From the definition of $\Omega_N^m$ and (\ref{eq estimate cm H1 H2}), we deduce that
\begin{eqnarray}\label{eq estimate nm L2}
\sup_{t\in[0,T]}\|n^m(t)\|^2_{L^2}+\int_0^T\|\nabla n^m(s)\|^2_{L^2}ds
\leq
C_{\mu,\delta,T,N,C_M,\|c(0)\|_{H^1}},\ \ \ \omega\in\Omega_N^m.
\end{eqnarray}
$\Box$
\vskip 0.3cm
\begin{Rem}\label{Rem 4.2}
(\ref{7-2}) and (\ref{eq estimate nm L2}) imply that
$$
\sup_{t\in[0,T]}\|n^m(t)\|^2_{L^2}+\int_0^T\|\nabla n^m(s)\|^2_{L^2}ds<\infty,\ \ \ P\text{-a.s..}
$$
\end{Rem}

\subsection{Existence of Martingale Weak Solutions}
\begin{definition}\label{def martingale solution}
We say that there exists a martingale weak solution to the system (\ref{eq system 00}) if
there exists
a stochastic basis $({\Omega},{\mathcal{F}},\{{\mathcal{F}}_t\}_{t\in[0,T]},{P})$ and, on this basis, a $U$- cylindrical
Wiener process ${W}$, a progressively measurable process $({n},{c},{u})$ satisfying
\begin{itemize}
\item[(1)] $P$-a.s.
\begin{eqnarray*}
 &&{n}(1+|x|+|\ln {n}|)\in L^\infty([0,T],L^1(\mathcal{O})),\ \ \nabla\sqrt{{n}}\in L^2([0,T],L^2(\mathcal{O})),\\
 && {c}\in L^{\infty}([0,T],L^\infty(\mathcal{O})\cap H^1(\mathcal{O}))\cap L^2([0,T],H^2(\mathcal{O})),\\
 && {u}\in C([0,T],H)\cap L^2([0,T],V);
\end{eqnarray*}

\item[(2)]
For all $\psi_1$, $\psi_2\in C^\infty([0,T]\times\mathcal{O})$ with compact supports with respect to the space variable, and $\psi_1(T,\cdot)=\psi_2(T,\cdot)=0$, $P$-a.s.
\begin{eqnarray*}
\int_{\mathcal{O}}\psi_1(0,x)n_0dx
=
\int_0^T\int_{\mathcal{O}}
     {n}[\partial_t\psi_1+\nabla\psi_1\cdot {u}+\delta\Delta\psi_1+\nabla\psi_1\cdot(\chi({c})\nabla {c})]
     dxdt,
\end{eqnarray*}
\begin{eqnarray*}
\int_{\mathcal{O}}\psi_2(0,x)c_0dx
=
\int_0^T\int_{\mathcal{O}}
    {c}[\partial_t\psi_2+\nabla\psi_2\cdot {u}+\mu\Delta\psi_2]-{n}k({c})\psi_2
    dxdt,
\end{eqnarray*}

\item[(3)]
For $e\in V$, $0\leq t\leq T$,
\begin{eqnarray*}
\langle {u}(t),e\rangle_{H,H}
&=&
\langle u_0,e\rangle_{H,H}-\nu\int_0^t \langle A {u}(s),e\rangle_{V^*,V}ds
-
\int_0^t \langle\mathcal{P}\{({u}(s)\cdot\nabla){u}(s)\},e\rangle_{V^*,V}ds \\
&&-\int_0^t \langle \mathcal{P}\{{n}(s)\nabla \phi\} ,e\rangle_{H,H}ds
 +
 \int_0^t\langle \sigma({u}(s))d{W}_s,e\rangle_{H,H}
\end{eqnarray*}
holds $P$-a.s..
\end{itemize}
\end{definition}
\begin{theorem}
Suppose the assumptions (A)-(C) in Section 2 hold. Then, there exists a martingale weak solution to the stochastic Chemotaxis-Navier-Stokes system (1.1).
\end{theorem}
{\bf Proof}.  Let $(n^m, c^m, u^m)$ be the solution constructed in Section 4.2. We will prove that the family $\{(n^m, c^m, u^m); m\geq 1\}$ is tight in the space
$L^2([0,T],L^2(\mathcal{O}))\times L^2([0,T],H^1(\mathcal{O}))\times L^2([0,T], H)$. To this end,  it suffices to show that the families $\{n^m; m\geq 1\}$, $\{c^m; m\geq 1\}$, $\{u^m; m\geq 1\}$  are respectively  tight in the spaces
$L^2([0,T],L^2(\mathcal{O}))$, $L^2([0,T],H^1(\mathcal{O}))$ and $L^2([0,T], H)$.

\vskip 0.3cm
Define
$$
\mathcal{Y}=\Big\{g\in L^2([0,T],H^2(\mathcal{O})),\ \frac{dg}{dt}\in L^2([0,T],L^2(\mathcal{O}))\Big\}
$$
and the norm
$$
\|g\|_{\mathcal{Y}}=\|g\|_{L^2([0,T],H^2(\mathcal{O}))}+\|\frac{dg}{dt}\|_{L^2([0,T],L^2(\mathcal{O}))}.
$$
By Chapter III Theorem 2.1 in \cite{Temam 1979}(Page 271) and the Kondrachov embedding theorem, it is known that
the embedding of $\mathcal{Y}$ into $L^2([0,T],H^1(\mathcal{O}))$ is compact.
By (\ref{eq estimate cm H1 H2}) and (\ref{7-5}), we get that
\begin{eqnarray*}
&&{P}\Big(\|c^m\|^2_{\mathcal{Y}}\leq 2 C_{\mu,T,N,C_M,\|c(0)\|_{H^1}}\Big)\\
&&\geq
{P}(\Omega_N^m)\\
&&\geq
1-\frac{3C}{N}.
\end{eqnarray*}
Since we can choose the integer $N$ as large as we wish, we conclude that the family $\{c^m\}$ is tight in $L^2([0,T],H^1(\mathcal{O}))$.
\vskip 0.3cm
Similarly,
(\ref{eq estimate nm dt}) and (\ref{eq estimate nm L2}) imply that $\{n^m\}$ is tight in $L^2([0,T],L^2(\mathcal{O}))$.

\vskip 0.3cm

Given $\kappa\in(0,1)$, let $W^{\kappa,2}([0,T],V^*)$ be the Sobolev space of all $g\in L^2([0,T],V^*)$ such that
$$
\int_0^T\int_0^T\frac{\|g(t)-g(s)\|^2_{V^*}}{|t-s|^{1+2\kappa}}dtds<\infty,
$$
endowed with the norm
$$
\|g\|^2_{W^{\kappa,2}([0,T],V^*)}:=\int_0^T\|g(t)\|^2_{V^*}dt+\int_0^T\int_0^T\frac{\|g(t)-g(s)\|^2_{V^*}}{|t-s|^{1+2\kappa}}dtds.
$$
Set
$$
B(u,u):=\mathcal{P}(u\cdot\nabla)u.
$$
Recall
$$
\langle B(u,u),v\rangle_{V^*,V}\leq C\|u\|_H\|u\|_V\|v\|_V,\ \ \ u,v\in V.
$$
It is known (see e.g.  \cite{visik-Fursikov}) that $B$ can be extended to a continuous operator
\begin{eqnarray}\label{B-ext}
B:H\times H\rightarrow D(A^{-\varrho})
\end{eqnarray}
for some $\varrho>1$.

Using the equation satisfied by $u^m$, applying the similar arguments as in the proof of Theorem 3.1 in \cite{Flandoli-Gatarek}, we can show that
\begin{eqnarray}\label{eq um tight L2 00}
\mathbb{E}\Big(\|u^m\|_{W^{\kappa,2}([0,T],V^*)}\Big)\leq C_\kappa.
\end{eqnarray}
Recall that the embedding of $L^2([0,T],V)\cap W^{\kappa,2}([0,T],V^*)$ into $L^2([0,T],H)$ is compact (see e.g. Theorem 2.1 in \cite{Flandoli-Gatarek}).  (\ref{eq um tight L2 00}) and (\ref{eq strong u}) imply that $\{u^m\}$
 is tight in $L^2([0,T],H)$.

 \vskip 0.3cm

 On the other hand, by (d) and (e) in Corollary \ref{cor 4.1}, applying Theorem 1 in \cite{D} and Corollary 5.2 in \cite{J}, as the proof of Lemma 4.5
 in \cite{Zhai Zhang}, we can prove that $\{u^m\}$ is tight in $D([0,T],D(A^{-\varrho}))$, here $\varrho$ is the constant appeared in (\ref{B-ext}) and $D([0,T],D(A^{-\varrho}))$ denotes the space of right continuous functions with left limits from $[0,T]$ into $D(A^{-\varrho})$
 equipped with the Skorokhod topology. Moreover, since $u^m$ takes values in $C([0,T],H)$,
Proposition 1.6 in \cite{J} implies that $\{u^m\}$ is also tight in $C([0,T],D(A^{-\varrho}))$ equipped with the usual uniform topology.
\vskip 0.3cm

Now we have proved that $\{(n^m,c^m,u^m)\}$ is tight in the space:
$$
\Pi:=L^2([0,T],L^2(\mathcal{O}))\times L^2([0,T],H^1(\mathcal{O}))\times C([0,T],D(A^{-\varrho}))\cap L^2([0,T],H).
$$
By the Skorohod embedding theorem, there exist a stochastic basis $(\widetilde{\Omega},\widetilde{\mathcal{F}}, \{\widetilde{\mathcal{F}}_t\}_{t\in[0,T]},\widetilde{\mathbb{P}})$ and, on this basis,
$\Pi$-valued random variables $(\tilde{n}^m,\tilde{c}^m,\tilde{u}^m)$, $(\tilde{n},\tilde{c},\tilde{u})$ such that
\begin{itemize}
\item[(I)] $(\tilde{n}^m,\tilde{c}^m,\tilde{u}^m)$ has the same law as $(n^m,c^m,u^m)$,

\item[(II)] $(\tilde{n}^m,\tilde{c}^m,\tilde{u}^m)\rightarrow (\tilde{n},\tilde{c},\tilde{u})$  in $\Pi$,  $\widetilde{\mathbb{P}}$-a.s.
\end{itemize}
By a similar argument as in the proof of Theorem 3.1 in \cite{Flandoli-Gatarek}, we can show that
$$
\tilde{u}(\cdot,\tilde{\omega})\in C([0,T],D(A^{-\varrho}))\cap L^\infty([0,T],H)\cap L^2([0,T],V),\ \ \widetilde{\mathbb{P}}\text{-a.s.},
$$
and there exists a $U$-cylindrical Wiener process
$\widetilde{W}$ on the stochastic basis $(\widetilde{\Omega},\widetilde{\mathcal{F}},\{\widetilde{\mathcal{F}}_t\}_{t\in[0,T]},\widetilde{\mathbb{P}})$ such that $\widetilde{\mathbb{P}}$-a.s., the identity
\begin{eqnarray}\label{eq tilde u 00}
\langle \tilde{u}(t),e\rangle_{H,H}
&=&
\langle u_0,e\rangle_{H,H}-\nu\int_0^t \langle A \tilde{u}(s),e\rangle_{V^*,V}ds
-
\int_0^t \langle\mathcal{P}\{(\tilde{u}(s)\cdot)\tilde{u}(s)\},e\rangle_{V^*,V}ds\nonumber \\
&&-
\int_0^t \langle\mathcal{P}\{\tilde{n}(s)\nabla \phi\},e\rangle_{H,H} ds
 +
 \int_0^t\langle \sigma(\tilde{u}(s))d\widetilde{W}_s,e\rangle_{H,H}
\end{eqnarray}
holds for all $t\in[0,T]$ and all $e\in D(A^\varrho)$. Furthermore, it follows  from (\ref{eq tilde u 00}) that $\widetilde{u}(\cdot)\in C([0,T],H)$ $\widetilde{\mathbb{P}}$-a.s..
This can be seen as follows.  Let $\widetilde{L}(t)$ be the solution of the stochastic evolution equation:
$$
d\widetilde{L}(t)=A\widetilde{L}(t)dt+\sigma(\tilde{u}(t))d\widetilde{W}_t,
$$
Then it is well known (see e.g. \cite{PR}) that $\widetilde{L}(\cdot,\tilde{\omega})\in C([0,T],H)\cap L^2([0,T],V)$ $\widetilde{\mathbb{P}}$-a.s..
On the other hand, using classical PDE arguments, Theorem 3.1 and Theorem 3.2 in \cite{Temam 1979}, we can show that there exists a unique process $\widetilde{Z}\in C([0,T],H)\cap L^2([0,T],V)$ satisfying the random PDE:
\begin{eqnarray}\label{eq widetile Z}
&& d\widetilde{Z}(t)+ \Big((\widetilde{Z}(t)+\widetilde{L}(t))\cdot\nabla\Big)(\widetilde{Z}(t)+\widetilde{L}(t))
 =
  \nu\Delta\widetilde{Z}(t)ds-\tilde{n}(t)\nabla \phi\\
 &&\widetilde{Z}(0)=u_0.\nonumber
\end{eqnarray}
From the equation (\ref{eq tilde u 00}), it is easy to see that $\widetilde{u}=\widetilde{L}+\widetilde{Z}$.  Hence
\begin{eqnarray}\label{eq tilde u in C}
\tilde{u}(\cdot,\tilde{\omega})\in C([0,T],H)\cap L^2([0,T],V),\ \ \widetilde{\mathbb{P}}\text{-a.s.}.
\end{eqnarray}
By a density argument,  it is easy to see that the identity (\ref{eq tilde u 00}) holds for all $e\in V$.

\vskip 0.3cm
Using the equations satisfied by $(n^m,c^m,u^m)$, we see that,
for all $\psi_1$, $\psi_2\in C^\infty([0,T]\times\mathcal{O})$ with compact supports with respect to the space variable, and $\psi_1(T,\cdot)=\psi_2(T,\cdot)=0$,
\begin{eqnarray*}
\int_{\mathcal{O}}\psi_1(0,x)n^m_0dx
=
\int_0^T\int_{\mathcal{O}}
     \tilde{n}^m[\partial_t\psi_1+\nabla\psi_1\cdot \tilde{u}^m+\delta\Delta\psi_1+\nabla\psi_1\cdot(\chi(\tilde{c}^m)\nabla \tilde{c}^m)]
     dxdt,
\end{eqnarray*}
\begin{eqnarray*}
\int_{\mathcal{O}}\psi_2(0,x)c^m_0dx
=
\int_0^T\int_{\mathcal{O}}
    \tilde{c}^m[\partial_t\psi_2+\nabla\psi_2\cdot \tilde{u}^m+\mu\Delta\psi_2]-\tilde{n}^mk(\tilde{c}^m)\psi_2
    dxdt.
\end{eqnarray*}
Taking $m$ into $\infty$ in the above two equations, we see that $(\tilde{n},\tilde{c},\tilde{u})$ satisfies (2) in Definition \ref{def martingale solution}.

\vskip 0.3cm
Finally, by (a) (b) (c) (d) in Corollary \ref{cor 4.1}, (I) (II), (\ref{eq estimate cm H1 H2}) and (\ref{eq tilde u in C}), we see that $(\tilde{n},\tilde{c},\tilde{u})$ satisfies (1) in Definition \ref{def martingale solution}. Hence, $(\tilde{n},\tilde{c},\tilde{u})$ is a martingale weak solution.  $\Box$
\subsection{Pathwise Weak Solution}
\begin{theorem}
Assume, in addition, that the function $\chi(\cdot)$ is a positive constant. Then there exists a unique  pathwise weak solution to  the stochastic Chemotaxis-Navier-Stokes system (1.1).
\end{theorem}
{\bf Proof}. From Theorem 4.1, we already know that there exists a martingale weak solution to system (1.1). By the Watanable and Yamada Theorem, we will complete the proof of the theorem if we can  show the pathwise uniqueness of the solutions. That is what we will do in the remaining part of the proof. Without loss of generality, we assume $\chi(\cdot)\equiv 1$.

\vskip 0.3cm

Assume that $(n_1,c_1,u_1)$ and $(n_2,c_2,u_2)$ are two solutions of the system (1.1) on the same  probability basis $(\Omega,\mathcal{F},\{\mathcal{F}_t\}_{t\in[0,T]},\mathbb{P})$, with a same $U$-valued cylindrical Wiener process $W$. We will prove that
$$
(n_1,c_1,u_1)=(n_2,c_2,u_2).
$$

For simplicity, set
$$
n^\Delta=n_1-n_2,\ \ \ c^\Delta=c_1-c_2,\ \ \ u^\Delta=u_1-u_2.
$$

By chain rule,
\begin{eqnarray}\label{eq martingale unique n 00}
&&  \|n^\Delta(t)\|^2_{L^2}+2\delta\int_0^t\|\nabla n^\Delta(s)\|^2_{L^2}ds\nonumber\\
&=&
    2\int_0^t\langle u_1(s)n_1(s)-u_2(s)n_2(s),\nabla n^\Delta(s)\rangle_{L^2,L^2}ds\nonumber\\
      &&+
    2\int_0^t\langle n_1(s)\nabla c_1(s)-n_2(s)\nabla c_2(s),\nabla n^\Delta(s)\rangle_{L^2,L^2}ds\nonumber\\
&\leq&
     \delta\int_0^t\|\nabla n^\Delta(s)\|^2_{L^2}ds
     +
     \frac{2}{\delta}\int_0^t\|u_1(s)n_1(s)-u_2(s)n_2(s)\|^2_{L^2}ds\nonumber\\
     &&+
     \frac{2}{\delta}\int_0^t\|n_1(s)\nabla c_1(s)-n_2(s)\nabla c_2(s)\|^2_{L^2}ds\nonumber\\
&=&\delta\int_0^t\|\nabla n^\Delta(s)\|^2_{L^2}ds
     +I^n_1(t)+I^n_2(t).
\end{eqnarray}
By (\ref{eq P2 01 um}) and (\ref{eq GNS L4}), for $I^n_1(t)$, we have for $\varepsilon >0$,
\begin{eqnarray}\label{eq MU n I1}
    I_1^n(t)
&\leq&
    C_\delta\Big(
              \int_0^t\|u_1(s)\|^2_{L^4}\| n^\Delta(s)\|^2_{L^4}ds
              +
              \int_0^t\| u^\Delta(s)\|^2_{L^4}\|n_2(s)\|^2_{L^4}ds
            \Big)\nonumber\\
&\leq&
   C_\delta\Big(
              \int_0^t\|u_1(s)\|_{L^2}\|\nabla u_1(s)\|_{L^2}
                                   \Big(
                                    \|\nabla n^\Delta(s)\|_{L^2}\| n^\Delta(s)\|_{L^2}\nonumber\\
                                    &&\ \ \ \ \ \ \ \ \ \ \ \ \ \ \ \ \ \ \ \ \ \ \ \  \ \ \ \ \ \ \ \ \ \ \ \ \ \ \ \ \ \ \ \ \ \  \ \ \ \ \ \ \ \ \ +
                                    \| n^\Delta(s)\|_{L^2}^2
                                    \Big)ds\nonumber\\
              &&\ \ \ \ +
              \int_0^t\| u^\Delta(s)\|_{L^2}\|\nabla u^\Delta(s)\|_{L^2}
                                        \Big(
                                        \|n_2(s)\|_{L^2}\|\nabla n_2(s)\|_{L^2}+\|n_2(s)\|^2_{L^2}
                                        \Big)
                                        ds
            \Big)\nonumber\\
&\leq&
      \frac{\delta}{4}\int_0^t\|\nabla n^\Delta(s)\|^2_{L^2}ds\nonumber\\
      &&+
      C_\delta\int_0^t\| n^\Delta(s)\|^2_{L^2}\Big(
                                                 \|u_1(s)\|^2_{L^2}\|\nabla u_1(s)\|^2_{L^2}+\|u_1(s)\|_{L^2}\|\nabla u_1(s)\|_{L^2}
                                                 \Big)
              ds\nonumber\\
      &&+
      \epsilon\int_0^t\|\nabla( u^\Delta(s))\|^2_{L^2}ds\\
      &&+
      C_{\delta,\epsilon}\int_0^t\| u^\Delta(s)\|^2_{L^2}\Big(
                                                             \|n_2(s)\|^2_{L^2}\|\nabla n_2(s)\|^2_{L^2}+\|n_2(s)\|^4_{L^2}
                                                           \Big)
      ds.\nonumber
\end{eqnarray}
By (\ref{eq GNS L4}) and (\ref{eq P2 02 cm}), for $\varepsilon>0$ we have
\begin{eqnarray}\label{eq MU n I2}
  I^n_2(t)
&\leq&
  C_\delta\Big(
              \int_0^t\|n_1(s)\nabla c^\Delta(s)\|^2_{L^2}ds
              +
              \int_0^t\| n^\Delta(s)\nabla c_2(s)\|^2_{L^2}ds
          \Big)\nonumber\\
&\leq&
  C_\delta\Big(
              \int_0^t\|n_1(s)\|^2_{L^4}\|\nabla c^\Delta(s)\|^2_{L^4}ds
              +
              \int_0^t\| n^\Delta(s)\|^2_{L^4}\|\nabla c_2(s)\|^2_{L^4}ds
          \Big)\nonumber\\
&\leq&
  C_\delta\int_0^t\Big(
                      \|n_1(s)\|_{L^2}\|\nabla n_1(s)\|_{L^2}+\|n_1(s)\|^2_{L^2}
                  \Big)\nonumber\\
                  &&\ \ \ \  \ \ \ \cdot\Big(
                  \| c^\Delta(s)\|_{H^2}\|\nabla c^\Delta(s)\|_{L^2}+\|\nabla c^\Delta(s)\|^2_{L^2}
                  \Big)ds\nonumber\\
  &&+
  C_\delta\int_0^t\Big(
                      \| n^\Delta(s)\|_{L^2}\|\nabla  n^\Delta(s)\|_{L^2}+\| n^\Delta(s)\|^2_{L^2}
                  \Big)\nonumber\\
                  &&\ \ \ \ \ \ \ \ \ \  \ \ \ \cdot\Big(
                  \|c_2(s)\|_{H^2}\|\nabla c_2(s)\|_{L^2}+\|\nabla c_2(s)\|^2_{L^2}
                  \Big)ds\nonumber\\
&\leq&
    \epsilon\int_0^t\| c^\Delta(s)\|^2_{H^2}ds\nonumber\\
    &&+
    C_{\delta,\epsilon}\int_0^t\|\nabla c^\Delta(s)\|^2_{L^2}\nonumber\\
                         &&\ \ \ \ \cdot
                         \Big(
                           \|n_1(s)\|^2_{L^2}\|\nabla n_1(s)\|^2_{L^2}+\|n_1(s)\|^4_{L^2}
                           +
                           \|n_1(s)\|_{L^2}\|\nabla n_1(s)\|_{L^2}+\|n_1(s)\|^2_{L^2}
                         \Big)
                         ds\nonumber\\
    &&+
      \frac{\delta}{4}\int_0^t\|\nabla n^\Delta(s)\|^2_{L^2}ds\nonumber\\
     &&+
     C_\delta\int_0^t\| n^\Delta(s)\|^2_{L^2}\\
                      &&\cdot
                            \Big(
                               \|c_2(s)\|^2_{H^2}\|\nabla c_2(s)\|^2_{L^2}+\|\nabla c_2(s)\|^4_{L^2}
                               +
                               \|c_2(s)\|_{H^2}\|\nabla c_2(s)\|_{L^2}+\|\nabla c_2(s)\|^2_{L^2}
                             \Big)
              ds.\nonumber
\end{eqnarray}
Combining (\ref{eq martingale unique n 00})--(\ref{eq MU n I2}), we have
\begin{eqnarray}\label{eq MU n}
&&  \|n^\Delta(t)\|^2_{L^2}+\frac{\delta}{2}\int_0^t\|\nabla n^\Delta(s)\|^2_{L^2}ds\nonumber\\
&\leq&
    C_\delta\int_0^t\| n^\Delta(s)\|^2_{L^2}\Big(
                                                 \|u_1(s)\|^2_{L^2}\|\nabla u_1(s)\|^2_{L^2}+\|u_1(s)\|_{L^2}\|\nabla u_1(s)\|_{L^2}
                                                 \Big)
              ds\nonumber\\
      &&+
      \epsilon\int_0^t\|\nabla( u^\Delta(s))\|^2_{L^2}ds\\
      &&+
      C_{\delta,\epsilon}\int_0^t\| u^\Delta(s)\|^2_{L^2}\Big(
                                                             \|n_2(s)\|^2_{L^2}\|\nabla n_2(s)\|^2_{L^2}+\|n_2(s)\|^4_{L^2}
                                                           \Big)
      ds\nonumber\\
      &&
      +
      \epsilon\int_0^t\| c^\Delta(s)\|^2_{H^2}ds\nonumber\\
    &&+
    C_{\delta,\epsilon}\int_0^t\|\nabla c^\Delta(s)\|^2_{L^2}\nonumber\\
                         &&\ \ \ \ \cdot
                         \Big(
                           \|n_1(s)\|^2_{L^2}\|\nabla n_1(s)\|^2_{L^2}+\|n_1(s)\|^4_{L^2}
                           +
                           \|n_1(s)\|_{L^2}\|\nabla n_1(s)\|_{L^2}+\|n_1(s)\|^2_{L^2}
                         \Big)
                         ds\nonumber\\
    &&+
    C_\delta\int_0^t\| n^\Delta(s)\|^2_{L^2}\nonumber\\
                      &&\cdot
                            \Big(
                               \|c_2(s)\|^2_{H^2}\|\nabla c_2(s)\|^2_{L^2}+\|\nabla c_2(s)\|^4_{L^2}
                               +
                               \|c_2(s)\|_{H^2}\|\nabla c_2(s)\|_{L^2}+\|\nabla c_2(s)\|^2_{L^2}
                             \Big)
              ds.\nonumber
\end{eqnarray}

\vskip 0.3cm

Now we estimate $\|c^\Delta\|^2_{H^1}$. By the chain rule, we have
\begin{eqnarray}\label{eq MU c 00}
&&    \|c^\Delta(t)\|^2_{L^2}+2\mu\int_0^t\|\nabla c^\Delta(s)\|^2_{L^2}ds\nonumber\\
&\leq&
      2\int_0^t\langle u_1(s)c_1(s)-u_2(s)c_2(s),\nabla c^\Delta(s)\rangle_{L^2,L^2}ds\nonumber\\
      &&-
      2\int_0^t\langle k(c_1(s))n_1(s)-k(c_2(s))n_2(s),c^\Delta(s)\rangle_{L^2,L^2}ds\nonumber\\
&\leq&
     \frac{\mu}{2}\int_0^t\|\nabla c^\Delta(s)\|^2_{L^2}ds
     +
     \frac{2}{\mu}\int_0^t\|u_1(s)c_1(s)-u_2(s)c_2(s)\|^2_{L^2}ds\nonumber\\
     &&+
     \int_0^t\|k(c_1(s))n_1(s)-k(c_2(s))n_2(s)\|^2_{L^2}ds
     +
     \int_0^t\|c^\Delta(s)\|^2_{L^2}ds,
\end{eqnarray}
and
\begin{eqnarray}\label{eq MU nabla c 00}
&&  \|\nabla c^\Delta(t)\|^2_{L^2} + 2\mu \int_0^t\|\Delta c^\Delta(s)\|^2_{L^2}ds\nonumber\\
&=&
     -2\int_0^t\langle u_1(s)\cdot \nabla c_1(s)-u_2(s)\cdot \nabla c_2(s),\Delta c^\Delta(s)\rangle_{L^2,L^2}ds\nonumber\\
     &&-
     2\int_0^t\langle k(c_1(s))n_1(s)-k(c_2(s))n_2(s),\Delta c^\Delta(s)\rangle_{L^2,L^2}ds\nonumber\\
&\leq&
     \mu \int_0^t\|\Delta c^\Delta(s)\|^2_{L^2}ds
     +
     \frac{2}{\mu}\int_0^t\|u_1(s)\cdot \nabla c_1(s)-u_2(s)\cdot \nabla c_2(s)\|^2_{L^2}ds\nonumber\\
     &&+
     \frac{2}{\mu}\int_0^t\|k(c_1(s))n_1(s)-k(c_2(s))n_2(s)\|^2_{L^2}ds.
\end{eqnarray}
By the similar arguments as in the proof of (\ref{eq MU n I1}), we have
\begin{eqnarray}\label{eq MU c I1}
    &&\frac{2}{\mu}\int_0^t\|u_1(s)c_1(s)-u_2(s)c_2(s)\|^2_{L^2}ds\nonumber\\
&\leq&
      \frac{\mu}{2}\int_0^t\|\nabla c^\Delta(s)\|^2_{L^2}ds\nonumber\\
      &&+
      C_\mu\int_0^t\|c^\Delta(s)\|^2_{L^2}\Big(
                                                 \|u_1(s)\|^2_{L^2}\|\nabla u_1(s)\|^2_{L^2}+\|u_1(s)\|_{L^2}\|\nabla u_1(s)\|_{L^2}
                                                 \Big)
              ds\nonumber\\
      &&+
      \epsilon\int_0^t\|\nabla u^\Delta(s)\|^2_{L^2}ds\\
      &&+
      C_{\mu,\epsilon}\int_0^t\|u^\Delta(s)\|^2_{L^2}\Big(
                                                             \|c_2(s)\|^2_{L^2}\|\nabla c_2(s)\|^2_{L^2}+\|c_2(s)\|^4_{L^2}
                                                           \Big)
      ds.\nonumber
\end{eqnarray}
Furthermore, for $\varepsilon >0$,
\begin{eqnarray}\label{eq MU c I2}
&&     \int_0^t\|k(c_1(s))n_1(s)-k(c_2(s))n_2(s)\|^2_{L^2}ds\nonumber\\
&\leq&
       2\Big(
              \int_0^t\|(k(c_1(s))-k(c_2(s)))n_1(s)\|^2_{L^2}ds
              +
              \int_0^t\|k(c_2(s))n^\Delta(s)\|^2_{L^2}ds
             \Big) \nonumber\\
&\leq&
    \sup_{r\in[0,C_M]}|k'(r)|^2\int_0^t\||c^\Delta(s)|\cdot|n_1(s)|\|^2_{L^2}ds
    +
    \sup_{r\in[0,C_M]}|k(r)|^2\int_0^t \|n^\Delta(s)\|^2_{L^2}ds\nonumber\\
&\leq&
    C\int_0^t\|c^\Delta(s)\|^2_{L^4}\|n_1(s)\|^2_{L^4}ds
    +
    C \int_0^t \|n^\Delta(s)\|^2_{L^2}ds\nonumber\\
&\leq&
    C\int_0^t
                 \Big(
                 \|\nabla c^\Delta(s)\|_{L^2}\|c^\Delta(s)\|_{L^2}+\|c^\Delta(s)\|^2_{L^2}
                 \Big)
                 \Big(
                 \|\nabla n_1(s)\|_{L^2}\|n_1(s)\|_{L^2}+\|n_1(s)\|^2_{L^2}
                 \Big)
                 ds\nonumber\\
    &&+
    C \int_0^t \|n^\Delta(s)\|^2_{L^2}ds\nonumber\\
&\leq&
   \epsilon\int_0^t\|\nabla c^\Delta(s)\|^2_{L^2}ds
        +
    C \int_0^t \|n^\Delta(s)\|^2_{L^2}ds\\
      &&+
      C_\epsilon \int_0^t\|c^\Delta(s)\|^2_{L^2}
                       \Big(
                          \|\nabla n_1(s)\|^2_{L^2}\|n_1(s)\|^2_{L^2}+\|n_1(s)\|^4_{L^2}
                          +
                          \|\nabla n_1(s)\|_{L^2}\|n_1(s)\|_{L^2}+\|n_1(s)\|^2_{L^2}
                       \Big)
      ds,\nonumber
\end{eqnarray}
and


\begin{eqnarray}\label{eq MU nabla c I2}
&&    \frac{2}{\mu}\int_0^t\|u_1(s)\cdot \nabla c_1(s)-u_2(s)\cdot \nabla c_2(s)\|^2_{L^2}ds\nonumber\\
&\leq&
      C_\mu\int_0^t\|u^\Delta(s)\cdot\nabla c_1(s)\|^2_{L^2}ds
      +
      C_\mu\int_0^t\|u_2(s)\cdot\nabla c^\Delta(s)\|^2_{L^2}ds\nonumber\\
&\leq&
      C_\mu\int_0^t\|u^\Delta(s)\|^2_{L^4}\|\nabla c_1(s)\|^2_{L^4}ds
      +
      C_\mu\int_0^t\|u_2(s)\|^2_{L^4}\|\nabla c^\Delta(s)\|^2_{L^4}ds\nonumber\\
&\leq&
     C_\mu\int_0^t\|u^\Delta(s)\|_{L^2}\|\nabla u^\Delta(s)\|_{L^2}
                   \Big(
                       \|c_1(s)\|_{H^2}\|\nabla c_1(s)\|_{L^2}+\|\nabla c_1(s)\|^2_{L^2}
                   \Big)
                   ds\nonumber\\
      &&+
      C_\mu\int_0^t\|u_2(s)\|_{L^2}\|\nabla u_2(s)\|_{L^2}
                    \Big(
                       \|c^\Delta(s)\|_{H^2}\|\nabla c^\Delta(s)\|_{L^2}+\|\nabla c^\Delta(s)\|^2_{L^2}
                   \Big)
                   ds\nonumber\\
&\leq&
      \epsilon\int_0^t\|\nabla u^\Delta(s)\|^2_{L^2}ds\nonumber\\
      &&+
      C_{\epsilon,\mu}\int_0^t\|u^\Delta(s)\|^2_{L^2}
                   \Big(
                       \|c_1(s)\|^2_{H^2}\|\nabla c_1(s)\|^2_{L^2}+\|\nabla c_1(s)\|^4_{L^2}
                   \Big)
                   ds\nonumber\\
      &&+
      \epsilon\int_0^t\|c^\Delta(s)\|^2_{H^2}ds\nonumber\\
      &&+
      C_{\epsilon,\mu}\int_0^t\|\nabla c^\Delta(s)\|^2_{L^2}
                               \Big(
                                   \|u_2(s)\|^2_{L^2}\|\nabla u_2(s)\|^2_{L^2}+\|u_2(s)\|_{L^2}\|\nabla u_2(s)\|_{L^2}
                               \Big)
                               ds.
\end{eqnarray}
By (\ref{eq P2 03 estimate H2}), we have
\begin{eqnarray}\label{eq MU nabla c I3}
\epsilon\int_0^t\|c^\Delta(s)\|^2_{H^2}ds
\leq
\epsilon C\Big(\int_0^t\|c^\Delta(s)\|^2_{L^2}ds+\int_0^t\|\nabla c^\Delta(s)\|^2_{L^2}ds+\int_0^t\|\Delta c^\Delta(s)\|^2_{L^2}ds\Big).
\end{eqnarray}

Combining (\ref{eq MU c 00})--(\ref{eq MU nabla c I3}), we arrive at
\begin{eqnarray}\label{eq Mu c H1}
&&  \|c^\Delta(t)\|^2_{L^2}+\|\nabla c^\Delta(t)\|^2_{L^2}+(\mu-\epsilon(1+\frac{2}{\mu})-\epsilon C)\int_0^t\|\nabla c^\Delta(s)\|^2_{L^2}ds
      +
    (\mu-\epsilon C)\int_0^t\|\Delta c^\Delta(s)\|^2_{L^2}ds\nonumber\\
&\leq&
   2\epsilon\int_0^t\|\nabla u^\Delta(s)\|^2_{L^2}ds
   +
   C_\mu\int_0^t\|n^\Delta(s)\|^2_{L^2}ds
   \nonumber\\
   &&+
   C_{\mu,\epsilon}\int_0^t\|c^\Delta(s)\|^2_{L^2}
                         \Big(
                             1+\|u_1(s)\|^2_{L^2}\|\nabla u_1(s)\|^2_{L^2}+\|u_1(s)\|_{L^2}\|\nabla u_1(s)\|_{L^2}\nonumber\\
                             &&\ \ \ \ \ \ \ \ \ \ \ \ \ \ \ \ \ \ \ \ \ \ \ \ \ \ \ \ \ \ \ \ \ \ \ \ \ \ \ \ \ \ \ \ \ \ \ \ +
                             \|n_1(s)\|^2_{L^2}\|\nabla n_1(s)\|^2_{L^2}+\|n_1(s)\|^4_{L^2}\nonumber\\
                             &&\ \ \ \ \ \ \ \ \ \ \ \ \ \ \ \ \ \ \ \ \ \ \ \ \ \ \ \ \ \ \ \ \ \ \ \ \ \ \ \ \ \ \ \ \ \ \ \ +
                             \|n_1(s)\|_{L^2}\|\nabla n_1(s)\|_{L^2}+\|n_1(s)\|^2_{L^2}
                         \Big)
                         ds\nonumber\\
   &&+
   C_{\epsilon,\mu}\int_0^t\|u^\Delta(s)\|^2_{L^2}
                                  \Big(
                                     \sum_{i=1}^2
                                         \Big(
                                              \|c_i(s)\|^2_{L^2}\|\nabla c_i\|^2_{L^2}+\|c_i(s)\|^4_{L^2}
                                         \Big)
                                  \Big)
                                  ds\nonumber\\
   &&+
   C_{\epsilon,\mu}\int_0^t\|\nabla c^\Delta(s)\|^2_{L^2}
                            \Big(
                                 \|u_2(s)\|^2_{L^2}\|\nabla u_2(s)\|^2_{L^2}+\|u_2(s)\|_{L^2}\|\nabla u_2(s)\|_{L^2}
                            \Big)
                            ds.
\end{eqnarray}

Since
$$
2\Big|\langle \mathcal{P}\{(u_1(t)\cdot\nabla)u_1(t)-(u_2(t)\cdot\nabla)u_2(t)\},u^\Delta(t)\rangle_{V',V}\Big|
\leq
\nu \|\nabla u^\Delta(t)\|^2_{L^2}+C_\nu\|u_2(t)\|^4_{L^4}\|u^\Delta(t)\|^2_{L^2},
$$
by $\rm It\hat{o}$'s formula, we have
\begin{eqnarray}\label{eq Mu u 00}
&&  \|u^\Delta(t)\|^2_{L^2}+2\nu\int_0^t\|\nabla u^\Delta(s)\|^2_{L^2}ds\nonumber\\
&=&
    -2\int_0^t\langle\mathcal{P}\{ (u_1(s)\cdot\nabla)u_1(s)-(u_2(s)\cdot\nabla)u_2(s)\},u^\Delta(s)\rangle_{V',V}ds\nonumber\\
    &&-
    2\int_0^t\langle \mathcal{P}\{n_1(s)\nabla\phi-n_2(s)\nabla\phi\},u^\Delta(s)\rangle_{L^2,L^2}ds\nonumber\\
    &&+
    2\int_0^t\langle \sigma(u_1(s))-\sigma(u_2(s)),u^\Delta(s)\rangle_{L^2,L^2}dW(s)
    +
    2\int_0^t\|\sigma(u_1(s))-\sigma(u_2(s))\|^2_{\mathcal{L}^2_0}ds\nonumber\\
&\leq&
     \nu \int_0^t\|\nabla u^\Delta(s)\|^2_{L^2}ds+C_\nu\int_0^t\|u_2(s)\|^4_{L^4}\|u^\Delta(s)\|^2_{L^2}ds\nonumber\\
     &&+
     \int_0^t\|u^\Delta(s)\|^2_{L^2}ds+\|\nabla\phi\|_\infty^2\int_0^t\|n^\Delta(s)\|^2_{L^2}ds\nonumber\\
     &&+
     2\int_0^t\langle \sigma(u_1(s))-\sigma(u_2(s)),u^\Delta(s)\rangle_{L^2,L^2}dW(s)
    +
    C\int_0^t\|u^\Delta(s)\|^2_{L^2}ds\nonumber\\
&\leq&
   \nu \int_0^t\|\nabla u^\Delta(s)\|^2_{L^2}ds+C_\nu\int_0^t\|u_2(s)\|^2_{L^2}\|\nabla u_2(s)\|^2_{L^2}\|u^\Delta(s)\|^2_{L^2}ds\nonumber\\
     &&+
     C\int_0^t\|u^\Delta(s)\|^2_{L^2}ds+C_\phi\int_0^t\|n^\Delta(s)\|^2_{L^2}ds\nonumber\\
     &&+
     2\int_0^t\langle \sigma(u_1(s))-\sigma(u_2(s)),u^\Delta(s)\rangle_{L^2,L^2}dW(s).
\end{eqnarray}
\vskip 0.3cm

\noindent Set
\begin{eqnarray*}
\Lambda(s):=\|n^\Delta(s)\|^2_{L^2}+\|c^\Delta(s)\|^2_{L^2}+\|\nabla c^\Delta(s)\|^2_{L^2}+\|u^\Delta(s)\|^2_{L^2}.
\end{eqnarray*}
Choosing  $\epsilon$ sufficiently small, by (\ref{eq MU n})(\ref{eq MU nabla c I3}) (\ref{eq Mu c H1}) and (\ref{eq Mu u 00}), we get that
\begin{eqnarray}\label{eq MU F00}
    \Lambda(t)
\leq
    C_{\epsilon,\mu,\nu,\delta}\int_0^t \Lambda(s)\Xi(s)ds
    +
    2\int_0^t\langle \sigma(u_1(s))-\sigma(u_2(s)),u^\Delta(s)\rangle_{L^2,L^2}dW_s.
\end{eqnarray}
here
\begin{eqnarray*}
  \Xi(s)
&=&
  1
  +
  \sum_{i=1}^2\sum_{j=1}^2
              \Big(
                \|n_i(s)\|^j_{L^2}\|\nabla n_i(s)\|^j_{L^2}+\|n_i(s)\|^{2j}_{L^2}
                 +
                \|c_i(s)\|^j_{L^2}\|\nabla c_i(s)\|^j_{L^2}
                +\|c_i(s)\|^{2j}_{L^2}
                \nonumber\\
                 &&\ \ \ \ \ \ \ \ \ \ \ \ \ \ \ \ \ +
                 \|c_i(s)\|^j_{H^2}\|\nabla c_i(s)\|^j_{L^2}+\|\nabla c_i(s)\|^{2j}_{L^2}
                 +
                 \|u_i(s)\|^j_{L^2}\|\nabla u_i(s)\|^j_{L^2}
               \Big).
\end{eqnarray*}
By the Gronwall' lemma, we arrive at
\begin{eqnarray}\label{eq MU F01}
    \Lambda(t)
\leq
    2\sup_{s\in[0,t]}\Big|\int_0^s\langle \sigma(u_1(s))-\sigma(u_2(s)),u^\Delta(s)\rangle_{L^2,L^2}dW(s)\Big|\cdot e^{C_{\epsilon,\mu,\nu,\delta}\int_0^t \Xi(s)ds}
\end{eqnarray}
Define
\begin{eqnarray*}
\tau_i^N&=&\inf_{t\geq 0}\{
                       \sup_{s\in[0,t]}\|n_i(s)\|^2_{L^2} \bigvee \int_0^t\|\nabla n_i(s)\|^2_{L^2}ds
                       \bigvee
                       \sup_{s\in[0,t]}\|c_i(s)\|^2_{H^1} \bigvee \int_0^t\|c_i(s)\|^2_{H^2}ds\\
                       &&\ \ \ \ \ \ \ \ \ \ \ \ \ \bigvee
                       \sup_{s\in[0,t]}\|u_i(s)\|^2_{L^2} \bigvee \int_0^t\|\nabla u_i(s)\|^2_{L^2}ds
                       \geq N
                     \}\bigwedge T,\ \ \ i=1,2.
\end{eqnarray*}
Put $\tau_N=\tau^N_1\bigwedge\tau^N_2$. Because $(n_i,c_i,u_i),\ i=1,2$ satisfy (1) in Definition \ref{def martingale solution}, we see that
\begin{eqnarray}\label{eq subsec 4.4 01}
\sup_{s\in[0,T]}\|c_i(s)\|^2_{H^1} + \int_0^T\|c_i(s)\|^2_{H^2}ds
+
\sup_{s\in[0,T]}\|u_i(s)\|^2_{L^2} + \int_0^T\|\nabla u_i(s)\|^2_{L^2}ds
<\infty,\ \ \ \mathbb{P}\text{-a.s..}
\end{eqnarray}
Repeating the arguments in Subsection \ref{subsec 4.2}, we can get the following result (see Remark \ref{Rem 4.2}):
\begin{eqnarray}\label{eq subsec 4.4 02}
\sup_{s\in[0,T]}\|n_i(s)\|^2_{L^2} + \int_0^T\|\nabla n_i(s)\|^2_{L^2}ds
<\infty,\ \ \ \mathbb{P}\text{-a.s..}
\end{eqnarray}
(\ref{eq subsec 4.4 01}) and (\ref{eq subsec 4.4 02}) imply that
\begin{eqnarray}\label{eq tau N}
\tau_N \nearrow T,\ \ \mathbb{P}\text{-a.s..}
\end{eqnarray}
Replace $t$ by  $t\wedge \tau_N$ in (\ref{eq MU F01}) to get
\begin{eqnarray}\label{eq MU F02}
    \Lambda(t\wedge \tau_N)
&\leq&
    2\sup_{s\in[0,t\wedge \tau_N]}\Big|\int_0^s\langle (\sigma(u_1(s))-\sigma(u_2(s)))dW_s,u^\Delta(s)\rangle_{L^2,L^2}\Big|\cdot e^{C_{\epsilon,\mu,\nu,\delta}\int_0^{t\wedge \tau_N} \Xi(s)ds}\nonumber\\
&\leq&
   C_{\epsilon,\mu,\nu,\delta,N}
   \sup_{s\in[0,t]}\Big|\int_0^{s\wedge \tau_N}\langle (\sigma(u_1(s))-\sigma(u_2(s)))dW_s,u^\Delta(s)\rangle_{L^2,L^2}\Big|.
\end{eqnarray}
By BDG inequality,
\begin{eqnarray*}
      \mathbb{E}\Big(\sup_{t\in[0,T]}\Lambda(t\wedge \tau_N)\Big)
&\leq&
      C_{\epsilon,\mu,\nu,\delta,N}
    \mathbb{E}\Big(
      \sup_{s\in[0,T]}\Big|\int_0^{t\wedge \tau_N}\langle (\sigma(u_1(s))-\sigma(u_2(s)))dW_s,u^\Delta(s)\rangle_{L^2,L^2}\Big|
          \Big)\\
&\leq&
      C_{\epsilon,\mu,\nu,\delta,N}
    \mathbb{E}\Big(
      \Big|\int_0^{T\wedge \tau_N}\|u^\Delta(s)\|^4_{L^2}ds\Big|^{1/2}
          \Big)\\
&\leq&
     \frac{1}{2}\mathbb{E}\Big(\sup_{t\in[0,T]}\|u^\Delta(t\wedge \tau_N)\|^2_{L^2}\Big)
     +
     C_{\epsilon,\mu,\nu,\delta,N}
    \mathbb{E}\Big(
      \int_0^{T\wedge \tau_N}\|u^\Delta(s)\|^2_{L^2}ds
          \Big)\\
&\leq&
     \frac{1}{2}\mathbb{E}\Big(\sup_{t\in[0,T]}\|u^\Delta(t\wedge \tau_N)\|^2_{L^2}\Big)
     +
     C_{\epsilon,\mu,\nu,\delta,N}
      \int_0^{T}\mathbb{E}\Big(\sup_{s\in[0,t]}\|u^\Delta(s\wedge \tau_N)\|^2_{L^2}dt\Big).
\end{eqnarray*}
By the Gronwall's lemma,
we obtain
$$
\mathbb{E}\Big(\sup_{t\in[0,T]}\Lambda(t\wedge \tau_N)\Big)=0.
$$
Let $N\rightarrow \infty$ to obtain
$$
\mathbb{E}\Big(\sup_{t\in[0,T]}\Lambda(t)\Big)=0,
$$
which implies the uniqueness. $\Box$

\vskip0.5cm {\small {\bf  Acknowledgements}\ \  This work is partly supported by National Natural Science Foundation of China (No.11671372, No.11431014, No.11401557)},  the Fundamental Research Funds
for the Central Universities (No. WK 3470000008), and Key Research Program of Frontier Sciences CAS(No. QYZDB-SSW-SYS009).

\end{document}